\documentclass[12pt]{iopart}

\usepackage{iopams}
\usepackage{graphicx}

\newcommand{\tfrac}{\textstyle\frac}
\def\operatorname#1{\mathop{\mathrm{#1}}\nolimits}
\newcommand{\mean}{\operatorname{mean}}

\newcommand{\A}{\mathbf A}

\newcommand{\M}{\mathbf M}
\newcommand{\N}{\mathbf N}
\newcommand{\R}{\mathbf R}
\newcommand{\T}{\mathbf T}
\newcommand{\Tss}{\mathbf{T}_{ss}}
\newcommand{\pp}{\mathbf{e}_1}
\newcommand{\qq}{\mathbf{e}_2}
\newcommand{\tpp}{\tilde{\mathbf{e}}_1}
\newcommand{\tqq}{\tilde{\mathbf{e}}_2}
\newcommand{\Ts}{\T_s}
\newcommand{\Xt}{\mathbf X_t}
\newcommand{\X}{\mathbf X}
\newcommand{\Xs}{\mathbf X_s}
\newcommand{\Xss}{\mathbf X_{ss}}
\newcommand{\Tt}{\mathbf T_t}

\newtheorem{theorem}{Theorem}

\begin{document}

\title{Vortex Filament Equation for a Regular Polygon}

\author{Francisco de la Hoz$^1$ and Luis Vega$^{2,3}$}

\address{$^1$ Department of Mathematics and Statistics and Operations Research, Faculty of Science and Technology, University of the Basque Country UPV/EHU, Barrio Sarriena S/N, 48940 Leioa, Spain}

\address{$^2$ Department of Mathematics, Faculty of Science and Technology, University of the Basque Country UPV-EHU, Barrio Sarriena S/N, 48940 Leioa, Spain}

\address{$^3$ BCAM-Basque Center for Applied Mathematics, Alameda Mazarredo 14, 48009 Bilbao, Spain}

\eads{\mailto{francisco.delahoz@ehu.es} and \mailto{luis.vega@ehu.es}}

\begin{abstract}

\noindent In this paper, we study the evolution of the vortex filament equation (VFE),
$$
\Xt = \Xs \wedge \Xss,
$$

\noindent with $\X(s, 0)$ being a regular planar polygon. Using algebraic techniques, supported by full numerical simulations, we give strong evidence that $\X(s, t)$ is also a polygon at any rational time; moreover, it can be fully characterized, up to a rigid movement, by a generalized quadratic Gau{\ss} sum.

We also study the fractal behavior of $\X(0, t)$, relating it with the so-called Riemann's non-differentiable function, that was proved by Jaffard to be a multifractal.

\end{abstract}


\submitto{\NL}

\maketitle

\section{Introduction}

Given a curve $\X_0 : \mathbb{R} \longrightarrow \mathbb{R}^3$, we consider the geometric
flow
\begin{equation}
\label{e:binormal}
\Xt = \kappa\,\mathbf b,
\end{equation}

\noindent where $\kappa$ is the curvature and $\mathbf b$ the binormal component of the
Frenet-Serret formulas
\begin{equation}
    \label{e:frenet-serret}
\pmatrix{
\T \cr \mathbf n \cr \mathbf b
}_s =
\pmatrix{
0 & \kappa & 0 \cr -\kappa & 0 & \tau \cr 0 & -\tau & 0
} \cdot
\pmatrix{
\T \cr \mathbf n \cr \mathbf b
}.
\end{equation}

\noindent The flow can be expressed as
\begin{equation}
\label{e:xt}\Xt = \Xs\wedge\Xss,
\end{equation}

\noindent where $\wedge$ is the usual cross-product, $t$ is the time, and $s$ is the arc-length parameter. It appeared for the first time in 1906 in the Ph. D. thesis of Da Rios \cite{darios} and was rederived in 1965 by Arms and Hama \cite{arms} as an approximation of the dynamics of a vortex filament under the Euler equations. This model is usually known as the vortex filament equation (VFE); we refer the reader to \cite{batchelor} and \cite{saffman} for an analysis and discussion of its limitations.

Equations \eref{e:binormal} and \eref{e:xt} are also known as the binormal equation and the localized induction approximation (LIA), respectively. Some of their explicit solutions are the line,
the circle, and the helix. Since the tangent vector $\T =
\X_s$ remains with constant length, we can assume that $\T\in\mathbb S^2$, for all $t$. Differentiating \eref{e:xt}, we get the so-called Schr\"odinger map equation onto the sphere:
\begin{equation}
\label{e:schmap}\Tt = \T\wedge\Tss,
\end{equation}

\noindent which is a particular case of the Landau-Lifshitz equation for ferromagnetism \cite{landau}. Equation \eref{e:schmap} can be rewritten in a more geometric way as
\begin{equation}
\T_t = \mathbf J\mathbf D_s\mathbf \T_s,
\end{equation}
where $\mathbf D_s$ is the covariant derivative, and $\mathbf J$ is
the complex structure of the sphere. Written in this way, \eref{e:schmap} can be generalized to more general definition domains and images, as in \cite{delahoz2007}, where the hyperbolic plane $\mathbb H^2$ was chosen as the target space.

A relevant step forward in the understanding of \eref{e:schmap} was given by Hasimoto in \cite{hasimoto}, where a transformation (see Section \ref{s:hasimoto}) that relates \eref{e:xt} and \eref{e:schmap} with the nonlinear Schr\"odinger (NLS) equation was introduced:
\begin{equation}
\label{e:schr1}
\psi_t = \rmi\psi_{ss} + \rmi\left(\frac{1}{2}(|\psi|^2 + A(t))\right)\psi, \quad A(t)\in\mathbb R.
\end{equation}

\noindent Equation \eref{e:xt} is time reversible, i.e., if $\X(s, t)$ is a solution, so is $\X(-s, -t)$. Bearing in mind this fact, an important property of \eref{e:xt} (and, hence, of \eref{e:schmap}) is that it has a one-parameter family of regular self-similar solutions that develop a corner-shaped singularity at finite time. This was proved in \cite{gutierrez} for the Euclidean case, and in \cite{delahoz2007} for the hyperbolic case. The self-similar solutions are written in the form
\begin{equation}
\X(s,t)=\frac{1}{\sqrt t }\mathbf G\left(\frac{s}{\sqrt t }\right).
\end{equation}

\noindent As was proved by Buttke in \cite{buttke87}, $\mathbf G$ is a smooth curve characterized by the geometric properties
\begin{equation}
\kappa=a, \qquad \tau= \frac s2,
\end{equation}

\noindent with $a$ being a constant. Therefore, the shape of $\mathbf G$ is the one of a hairpin: it looks like a segment of a circle of radius $1/a$ close to the origin, that becomes two helicoidal curves with a bigger pitch for bigger $|s|$, and that winds around two different lines \cite{DelahozGarciaCerveraVega09}.

Buttke also studied numerically these self-similar solutions. Later on, in \cite{DelahozGarciaCerveraVega09}, a careful numerical study of those solutions for both the Euclidean and the hyperbolic cases was done. On the one hand, the authors reproduced numerically in \cite{DelahozGarciaCerveraVega09} the formation of the corner-shaped singularity, and, on the other hand, they started with a corner-shaped initial datum, recovering numerically the self-similar solutions; in all cases, the correct choice of boundary conditions was shown to be vital. Furthermore, they gave numerical evidence that fractality phenomena appeared in \eref{e:schmap} (see Figure 2.6 in p. 1059 of \cite{DelahozGarciaCerveraVega09}), if, for instance, fixed boundary conditions were imposed on the tangent vector $\T$. Remark that the relationship between fractals and the Schr\"odinger map equation is not new; indeed, in \cite{Peskin94}, an aortic valve model was proposed, in order to study the apparent fractal character of the valve's fiber architecture. Instead of $\X(s, t)$, the authors wrote $\X(u, v)$, where $v$, which corresponds to our time, is such that the curves $v = \mbox{const}$ are the fibers. After imposing fixed boundary conditions at $u = \pm u_0$, they showed that the curves $u = \mbox{const}$ have a fractal character. Besides, the fractal dimension of those curves was calculated numerically in \cite{Stern94}.

It is well known that this fractal behavior already occurs at the linear level, as can be easily seen when solving numerically the free Schr\"odinger equation
\noindent
\begin{equation}
\label{e:linear}
\psi_t=\frac{\rmi}{4\pi}\psi_{ss},
\end{equation}

\noindent with an initial datum that is piecewise constant as the sign function in, say, $s\in[-\pi, \pi]$, and then extended by periodicity. The trigonometric series associated to these problems have been extensively studied and, as we will see later, the behavior at rational times modulo $2\pi$ is completely different to that at irrational ones. In \cite{BG}, Berry and Goldberg made a far-reaching connection of this behavior with the so called Talbot effect in optics. In \cite{BK} and in \cite{BMS}, the authors used \eref{e:linear} to model the Talbot effect and showed that at rational times the solution can be obtained as a finite overlapping of translates of the initial datum, while at irrational times the images have a fractal profile. Also in \cite{B1}, a conjecture was made about the dimensions of those fractals. As far as we know, the only available results on the nonlinear setting are the recently ones obtained by Erdo{\u g}an and Tzirakis \cite{ET} (see also \cite{Ol}) in the cubic nonlinear Schr\"odinger equation (i.e. equation \eref{e:schr1}, with $A(t)=0$), for the initial data given by piecewise continuous functions. This means that the data are ``almost'' in the Sobolev class $\mathcal H^{1/2}$, that is subcritical with respect to the one determined by the scaling of the equation (i.e. $\mathcal H^{-1/2}$). This typically implies that the nonlinear potential can be seen as an external force. Then, its contribution is computed using Duhamel's integral, which has a regularizing effect. As a consequence, the fractal behavior is due to the linear part of the solution. At this point, the available results for the linear theory are used (see \cite{Os} and \cite{KR}). More recently, Erdo{\u g}an and Tzirakis, in collaboration with Chousionis, have also studied in \cite{CET} the Schr\"odinger map equation, but still with subcritical initial conditions.

In this paper, we are interested in filaments with corners. This means to consider delta functions as the initial data for \eref{e:schr1} or piecewise continuous functions as the initial data for \eref{e:schmap}. Therefore, we work with data in critical spaces and the problem turns out to be much more involved. The case of an infinite curve with just one corner and that is otherwise smooth has been studied in the sequence of papers \cite{BV1,BV2,BV3} by Banica and Vega. One of the conclusions is that \eref{e:schr1} is in fact ill-posed for this type of initial conditions, so that \eref{e:binormal} is also ill-posed if it is understood in a classical sense. However, after introducing appropriate function spaces and using the appropriate definition of solution, it is proved that \eref{e:xt} and \eref{e:schmap} are well-posed for this type of initial condition. More concretely, the velocity at the point where the corner is located is determined by the self-similar solution that at time zero has a corner with the same angle. In fact, in \cite{gutierrez} it was proved that there exists just one self-similar solution with this property.

Even if the solutions of \eref{e:xt} for an initial datum with a corner are well understood, nothing had been done for more general initial data with several corners, in particular, polygons. Nevertheless, in two recent papers \cite{didier,Didier2} by Jerrard and Smets, they prove a global existence theorem that allows to consider such types of initial data. Moreover, they simulate numerically the evolution of the unit square at different times, suggesting that the solution could become again polygonal; indeed, at a certain time, the square seems to reappear, but with the axes turned by a $\pi/4$-angle with respect to the initial ones. In the following pages, we will show not only that these observations are correct, but that, given a regular planar polygon with $M$ sides as the initial datum, another polygon is obtained at times of the form $t_{pq}=(2\pi/M^2)(p/q)$, $\gcd(p, q) = 1$. In fact, we will give a complete description of this polygon; in particular, and except for $q = 2$, we will see that it is a skew polygon with $Mq$ sides (if $q$ odd), or $Mq/2$ sides (if $q$ even).

One could wander what, if any, is the connection of these solutions with real fluids. Of course, one cannot expect anything at the quantitative level, because  Da Rios' and Arms and Hama's approximations are based first in a truncation of the Biot-Savart integral and then in a renormalization of the time variable. However, using VFE as an approximation of the evolution of real filaments cannot be ruled out at the qualitative level, even with data as singular as the ones we consider. In fact, in \cite{DelahozGarciaCerveraVega09}, the similitude between the self-similar solutions and a flow traversing a triangular wing is pointed out (see Figure 1.1 in  p. 1051 of that paper).

In the concrete case of the solutions considered in this article, we could ask if the macroscopic effect of the turning axes mentioned above can be seen in the dynamics of real fluids. It appears that this is something quite well documented for noncircular jets. It was first observed when the nozzle has an elliptic shape and, later on, seen for shapes given by a rectangle, a square, and an equilateral triangle among others. The interest of considering nozzles with corners started in the mid 80's and has received considerable attention since then. We refer the reader to the survey \cite{GG}, where the papers related to the above observations can be found, together with a detailed analysis on how sensitive the axes-turning phenomenon is to the initial geometric conditions; and several numerical experiments are also reviewed. It has to be understood that, at the numerical level, the corners are always effectively rounded and, in real experiments, the corners are immediately smoothed out once the jet leaves the nozzle, due to viscosity effects. This is consistent with the behavior of the self-similar solutions mentioned above, because the initial corner evolves into a smooth curve with the shape of a hairpin.

Furthermore, the appearance of smooth polygonal structures with more angles has  been appreciated for nozzles with the shape of an equilateral triangle (see \cite[Figure 6]{GG}). Also in \cite[Figure 10]{GGP}, while in \cite[p. 1492]{GGP} it is said  ``...a consistent eightfold distribution pattern is also suggested...'' for square nozzles. Similar results are obtained in the numerical simulations of \cite{MMG} (see Figure 8 in p. 9, and p. 11, where the $\pi$-degree ``flip-flop'' in the case of an equilateral triangle is emphasized). We think it would be worth studying in more detail whether or not these more complicated polygonal structures appear in noncircular jet flows.

The structure of this paper is as follows. In Section \ref{s:hasimoto}, we review the Hasimoto transformation for \eref{e:schmap}. In Section \ref{s:Xt}, we gather the main theoretical arguments that support our numerical experiments in the form of a theorem. We start recalling in Section \ref{s:formulation} some elementary geometric symmetries that regular polygons have and that are preserved by VFE. On the one hand, those symmetries greatly simplify the numerical implementation and, on the other hand, they play a fundamental role in computing some important quantities that govern the dynamics of regular polygons. Then, in Section \ref{s:galilean}, we pay especial attention to a group of symmetries that, as far as we know, can only be visualized through the use of the Hasimoto transformation and the NLS equation \eref{e:schr1}. They are the so-called Galilean transformations \eref{e:galilean}, and leave invariant the set of solutions of VFE.

We proceed as in \cite{KPV}, where the solution of the NLS equation that has the delta function as an initial condition is characterized by the fact that it leaves invariant all the Galilean transformations; this property determines the solution, except for one function of time that is easily computed by solving a simple ODE and choosing $A(t)$ in \eref{e:schr1} appropriately. In our case, we take a planar regular polygon with $M$ sides as the initial datum.  Then, as in \cite{KPV}, the solution of \eref{e:schr1} is determined except, again, by a function that depends on time; but, in this case, it is very delicate how to specify this function and does not seem to be easy. Our approach in this paper is to integrate the Frenet-Serret frame with
\begin{equation}
\label{e:psist0}
\psi(s,t)=\hat \psi(0,t)\sum_{k = -\infty}^{\infty} \rme^{-\rmi(Mk)^2t+\rmi Mks},
\end{equation}

\noindent where $\psi(s,t)$, which is in this case a distribution, is obtained from the Hasimoto transformation.

In Section \ref{s:rational}, we evaluate $\psi(s, t)$ at times $t_{pq}$ that are rational multiples of $2\pi/M^2$, i.e., $t_{pq}=(2\pi/M^2)(p/q)$, $\gcd(p,q)=1$, because $\psi(s, t_{pq})$ takes then the form \eref{e:psistpq}, so the corresponding curve is a skew polygon with $Mq$ sides, for $q$ odd, and with $Mq/2$ sides, for $q$ even. Remark that the structure of the polygon is determined by the generalized quadratic Gau{\ss} sums $G(-p, m, q)$ appearing in \eref{e:psistpq} (see \ref{s:gausssums}); in particular, when $q$ is even, half of the sums are zero, which ultimately explains why the polygon has $Mq/2$ sides in that case, instead of $Mq$.

The value $\hat \psi(0,t_{pq})$ is obtained in Section \ref{s:recovering}, by imposing the condition that the polygon has to be closed. As a consequence, the tangent vector $\T$ is determined, except for a rotation that can depend on time. This rotation is computed in Section \ref{s:rigid} by using the extra symmetries that are available for regular polygons and that we have already mentioned. The next step is to integrate $\T$, in order to obtain the corresponding skew polygon $\X$ at rational times. Using again the symmetries, we are able to calculate all the necessary quantities, except for some vertical translation of $\X$ that can be easily computed numerically (although we do not see a simple way of doing it by using just theoretical arguments).

In Section \ref{s:numerical}, we propose a numerical method for \eref{e:xt} and \eref{e:schmap}. In Section \ref{s:symmetries}, we explain carefully how to take advantage of the symmetries of $\X$ and $\T$ for the types of solutions considered, which leads to a dramatically improvement in the computational cost, as explained in Section \ref{s:cost}.

In Section \ref{s:experiments}, we perform numerical experiments, simulating the evolution of \eref{e:xt} and \eref{e:schmap} for different numbers of initial sides $M$. Concerning the vertical translation of $\X$, we give strong evidence that the center of mass of $\X$ moves upward with a constant speed $c_M\in(0, 1)$ that depends on $M$. We compare the numerical results with the algebraic values from Section \ref{s:Xt}, obtaining a complete agreement between both (totally different) approaches.

In Section \ref{s:riemann}, we study the evolution of $\X(s, t)$, for $s = 0$. At this point, it is important to recall the invariances of the so-called Jacobi theta function, closely related to \eref{e:psist0}:
\begin{equation}
\label{e:jacobi}\theta(s,t)=\sum_{k = -\infty}^{\infty}\rme^{-\pi \rmi k^2t+2\pi \rmi ks}.
\end{equation}

\noindent This function is precisely the solution of the free Schr\"odinger equation \eref{e:linear} with a one-periodic delta as the initial condition. Therefore, the Galilean transformations still hold in this case; moreover, they also imply the extra symmetry
\begin{equation}
\theta(s,t)=\rme^{\rmi\pi/4}\frac {\rme^{\frac{\pi \rmi s^2}{t}}}{\sqrt t}\theta\left(\frac{s}{t}, \frac{-1}{t}\right), \quad t>0,
\end{equation}

\noindent which, together with the periodicity property $\theta(s,t+2)=\theta(s,t)$, generates the so-called unimodular subgroup. In fact, in \cite{Du,Ja}, all these symmetries are used to prove the existence of a continuum range of exponents for the H\"older regularity of
\begin{equation}
\label{e:phit0}
\phi(t) = \sum_{k = 1}^{\infty}\frac{\rme^{\pi \rmi k^2 t}}{\rmi \pi k^2}, \quad t\in[0,2],
\end{equation}

\noindent whose real part is precisely Riemann's non-differentiable function. More concretely, in \cite{Ja}, it was proved that Riemann's non-differentiable function is a multifractal (see also \cite{CC}, where other relevant trigonometric sums are considered). In other words, it was proved that the set of times $t$ that have the same H\"older exponent is a fractal with a dimension depending on the H\"older exponent, in such a way that the conjecture stated by Frisch and Parisi in \cite{FP} is fulfilled (see also \cite {F}, for more details at this respect and the connection of this question with fully developed turbulence and intermittency).

We could  ask if something similar to the properties of Riemann's non-differentiable function also holds in our case. It is very easy to find the analogous formulas to \eref{e:phit0} in out setting. In fact, for a given number of initial sides $M$, if we write $\X(0, t) = (X_1(0,t), X_2(0,t), X_3(0,t))$, bearing in mind the symmetries of the problem, we conclude that $\X(0, t)$ is a planar curve, so we identify the plane where it lives with $\mathbb C$ and define
\begin{equation}
\label{e:z(t)}
z(t) = -\|(X_1(0, t), X_2(0, t))\| + iX_3(0, t).
\end{equation}

\noindent This curve, or, more precisely, $z(t) - c_Mt$, can be seen as a nonlinear version of \eref{e:phit0}. Is $z(t)$ a multifractal? If the answer is positive, what is its spectrum of singularities in the sense of \cite{Ja}? We deem these two questions rather challenging from both a numerical and an analytical point of view. The results of Section \ref{s:riemann} strongly suggest that the answer to the first question is positive. In fact, after computing $z(t) - c_Mt$ for different $M$, and measuring the difference in the $\mathcal L^\infty$-norm between it and an appropriately scaled and rotated version of $\phi(t)$, the convergence rate is rather strong. We also give some numerical evidence that the H\"older exponent of $z(t)$ is $1/2$ for rational times $t_{pq}$; nevertheless, how the constants depend of the denominator $q$, which is a fundamental ingredient in the arguments in \cite{Du,Ja}, is unclear. This question deserves a much more detailed analysis that we plan to make in a forthcoming paper.

In Section \ref{s:fractalsT}, we study some interesting fractality phenomena associated to the tangent vector $\T$. In Section \ref{s:conclusions}, we draw the main conclusions, as well as some open questions that we postpone for the future. Finally, to conclude this paper, we offer in \ref{s:gausssums} a detailed study on the generalized quadratic Gau{\ss} sums.

\section{The Hasimoto transformation}

\label{s:hasimoto}

A central point of this paper is the natural connection between \eref{e:xt}-\eref{e:schmap} and the NLS equation. Indeed, applying the Hasimoto transformation \cite{hasimoto} to \eref{e:frenet-serret}:
\begin{equation}
\psi(s, t) = \kappa(s, t)\exp\left(\rmi\int^s\tau(s',t)ds'\right),
\end{equation}

\noindent the function $\psi$ satisfies the NLS equation
\begin{equation}
\label{e:schr}
\psi_t = \rmi\psi_{ss} + \rmi\left(\frac{1}{2}(|\psi|^2 + A(t))\right)\psi,
\end{equation}

\noindent where $A(t)$ is a certain real constant that depends on time. Nevertheless, for our purposes, it is not convenient to work with the torsion $\tau$. Instead, we consider another version of the Frenet-Serret trihedron; it is easy to check that all its possible generalizations have the form
\begin{equation}
\label{e:Te1e2}
\pmatrix{
\T \cr \pp \cr \qq
}_s =
\pmatrix{
    0 & \alpha & \beta
    \cr
    - \alpha & 0 & \gamma \cr -\beta & -\gamma
& 0
} \cdot
\pmatrix{
    \T \cr \pp \cr \qq
},
\end{equation}

\noindent for some vectors $\pp$ and $\qq$ that form an orthonormal base with the tangent vector $\T$. Moreover, we can make, without loss of generality, one of the coefficients $\alpha$, $\beta$, or $\gamma$ equal
zero. If we choose $\beta = 0$, denoting $\alpha \equiv \kappa$ and
$\gamma \equiv \tau$, we recover the habitual trihedron. On the other hand, in this paper, we choose $\gamma \equiv 0$, in order to avoid working with the torsion. In that case, the Hasimoto transformation takes the form
\begin{equation}
\psi \equiv \alpha + \rmi\beta,
\end{equation}

\noindent where
\begin{equation}
\label{e:ab}
\eqalign{
\alpha(s,t) = \kappa(s, t)\cos\left(\int^s\tau(s',t)ds'\right),
    \\
\beta(s,t) = \kappa(s, t)\sin\left(\int^s\tau(s',t)ds'\right).
}
\end{equation}

\noindent It is quite straightforward to check that this new definition of $\psi$ also satisfies \eref{e:schr}. The key point is to combine the basis vectors $\pp$ and $\qq$ into $\N \equiv \pp + \rmi\qq$. Since the proof closely matches that in \cite{hasimoto}, we omit here the details.

The main idea is to work with \eref{e:schr}, and, at a given $t$, recover the curve $\X(s, t)$ and the tangent vector $\T(s, t)$ from $\psi(s, t)$ by integrating \eref{e:Te1e2}, up to a rigid movement that we determine by the symmetries of the problem. Observe that, if we define
\begin{equation}
\label{e:hasimoto2}
\tilde\psi \equiv \rme^{\rmi\omega_0}\psi = \rme^{\rmi\omega_0}(\alpha + \rmi\beta), \quad \omega_0\in\mathbb R,
\end{equation}

\noindent and integrate \eref{e:Te1e2} for $\tilde\psi$, the new solution of \eref{e:Te1e2} is $\{\tilde\T, \tilde\pp, \tilde\qq\}$, where $\tilde \T \equiv \T$, and $\tilde\pp + \rmi\tilde\qq \equiv \rme^{\rmi\omega_0}(\pp + \rmi\qq)$, i.e., $\T$ does not change.

\section{A solution of $\Xt = \Xs\wedge\Xss$ for a regular polygon}

\label{s:Xt}

The aim of this paper is to understand the evolution of \eref{e:xt} for polygonal initial data. In fact, we limit ourselves to studying the simplest case, i.e., that of a regular planar polygon with $M$ sides. More precisely, we prove the following theorem:

\begin{theorem}

\label{t:theorem}

Assume that there is a unique solution of the initial value problem
\begin{equation}
\Xt = \Xs \wedge \Xss,
\end{equation}

\noindent where $\X(s,0)$ is a regular planar polygon with $M$ sides. Then, at a time $t_{pq}$ that is a rational multiple of $2\pi/M^2$, i.e., $t_{pq} \equiv (2\pi/M^2)(p/q)$, with $p\in\mathbb Z$, $q\in\mathbb N$, $\gcd(p, q) = 1$, the solution is a skew polygon with $Mq$ sides (if $q$ odd) or $Mq/2$ sides (if $q$ even). All the new sides have the same length; and the angle $\rho$ between two adjacent sides is constant. Furthermore, the polygon at a time $t_{pq}$ is the solution of the Frenet-Serret system
\begin{equation}
\label{e:Te1e2ab}
\pmatrix{
\T(s, t_{pq}) \cr \pp(s, t_{pq}) \cr \qq(s, t_{pq})
}_s =
\pmatrix{
    0 & \alpha(s, t_{pq}) & \beta(s, t_{pq})
    \cr
    - \alpha(s, t_{pq}) & 0 & 0 \cr -\beta(s, t_{pq}) & 0
& 0
} \cdot
\pmatrix{
    \T(s, t_{pq}) \cr \pp(s, t_{pq}) \cr \qq(s, t_{pq})
}
,
\end{equation}

\noindent where $\alpha(s, t_{pq}) + \rmi\beta(s, t_{pq}) = \psi(s, t_{pq})$, and $\psi(s, t_{pq})$ is the $2\pi/M$-periodic function defined over the first period $s\in[0, 2\pi/M)$ as
\begin{equation}
\label{e:psistpq0}
\psi(s, t_{pq}) =
\cases{
\frac{\rho}{\sqrt q}\sum_{m = 0}^{q - 1}G(-p, m, q) \delta(s - \tfrac{2\pi m}{Mq}), & if $q$ odd,
    \\
\frac{\rho}{\sqrt{2q}}\sum_{m = 0}^{q - 1}G(-p, m, q) \delta(s - \tfrac{2\pi m}{Mq}), & if $q$ even,
}
\end{equation}

\noindent with
\begin{equation}
\label{e:Gabc}
G(a, b, c) = \sum_{l=0}^{c - 1}\rme^{2\pi \rmi (al^2 + bl)/c},  \quad a, b\in\mathbb Z, c\in\mathbb Z \setminus \{0\},
\end{equation}

\noindent being a generalized quadratic Gau{\ss} sum (see \ref{s:gausssums}).

\end{theorem}

The rest of this section is devoted to the proof of this theorem. Moreover, we also conjecture that the angle $\rho$ is given by
\begin{equation}
\label{e:cosrho}
\cos(\rho) =
\cases{
2\cos^{2/q}(\tfrac{\pi}{M}) - 1, & if $q$ odd,
    \cr
2\cos^{4/q}(\tfrac{\pi}{M}) - 1, & if $q$ even.
}
\end{equation}

\noindent \textbf{Remark.} From our point of view, uniqueness is a very challenging problem for this type of equations, because of the critical regularity that we are assuming in the initial data. As we say in the introduction, uniqueness has been established in a series of papers by Banica and Vega, when the initial datum has just one corner and is otherwise regular. Banica and Vega prove that uniqueness holds for the binormal equation and for the Schr\"odinger map equation, but the cubic NLS equation is ill-posed as an initial value problem with initial datum being the Dirac delta function. Therefore, the work of Erdo{\u g}an and Tzirakis \cite{ET} cannot be used in our case, and the only general result for the binormal flow for data as singular as ours is so far the one by Jerrard and Smets \cite{Didier2}, where, on the other hand, they just prove existence.

\subsection{Formulation of the problem}

\label{s:formulation}

Since \eref{e:xt} and \eref{e:schmap} are invariant with respect to rotations, given a regular planar polygon with $M$ sides, we can assume without loss of generality that $\X(s, 0)$, and hence $\T(s, 0)$, live in the plane $OXY$. Identifying $OXY$ with $\mathbb C$, and assuming without loss of generality a total length of $2\pi$, $\X(s, 0)$ is the polygon parameterized by arc-length whose $M$ vertices, located at $s_k = 2\pi k / M$, $k = 0, \ldots, M-1$, are
\begin{equation}
\label{e:xt0}
\X(s_k, 0) = \frac{-\rmi \pi \rme^{\rmi \pi (2k - 1) / M}}{M\sin(\pi / M)},
\end{equation}

\noindent and a point $\X(s, 0)$, for $s_k < s < s_{k + 1}$, lays in the segment that joints $\X(s_k, 0)$ and $\X(s_{k + 1}, 0)$. Observe that it is straightforward to extend $\X(s, 0)$ periodically to the whole real with period $2\pi$. Then, since $\X(s_{k + 1}, 0) - \X(s_k, 0) = (2\pi/M)\rme^{2\pi \rmi k/M}$, the corresponding tangent vector is the periodic function with period $2\pi$ such that
\begin{equation}
\label{e:tt0}
\T(s, 0) = \rme^{2\pi \rmi k/M}, \quad \mbox{for } s_k < s < s_{k + 1}.
\end{equation}

\noindent The curve $\X(s, 0)$ is continuous, while its tangent vector $\T(s, 0)$ is only piecewise continuous. In fact, $\X(s, 0)$ can be regarded as a curve whose $2\pi/M$-periodic curvature is a sum of Dirac deltas:
\begin{equation}
\label{e:curvatureM}
\kappa(s) = \frac{2\pi}{M}\sum_{k = -\infty}^{\infty}\delta(s - \tfrac{2\pi k}{M});
\end{equation}

\noindent and the constant $2\pi/M$ has been chosen in order for the integral of the curvature over an interval of length $2\pi$ to be equal to $2\pi$.

At this point, it is necessary to understand the role played by symmetries in \eref{e:xt} and \eref{e:schmap} for these initial data. Both equations are invariant with respect to rotations, i.e, given a rotation matrix $\R$, if $\X = (X_1, X_2, X_3)^T$ and $\T = (T_1, T_2, T_3)^T$ are, respectively, solutions of them, so are $\R\cdot \X$ and $\R\cdot \T$ (in this paper, all the vectors are represented in column form). Therefore, if $\R\cdot \X(s, 0) = \X(s, 0)$ and $\R\cdot \T(s, 0) = \T(s, 0)$, then, if the solution is unique, $\R\cdot \X(s, t) = \X(s, t)$ and $\R\cdot \T(s, t) = \T(s, t)$, for all $t$. In particular, since $\X(s, 0)$ and $\T(s, 0)$, as defined in \eref{e:xt0} and \eref{e:tt0}, are invariant with respect to rotations of angle $2\pi k/ M$ around the $z$-axis, for all $k\in\mathbb Z$, we conclude that also $\X(s, t)$ and $\T(s, t)$ are invariant with respect to that kind of rotations, for all $t$. More precisely,
\begin{equation}
\label{e:rotationXT}
\eqalign{
X_1(s + \tfrac{2\pi k}{M}, t) + \rmi X_2(s + \tfrac{2\pi k}{M}, t) = \rme^{2\pi \rmi  k/M}(X_1(s, t) + \rmi X_2(s, t)),
    \cr
X_3(s + \tfrac{2\pi k}{M}, t) = X_3(s, t),
    \cr
T_1(s + \tfrac{2\pi k}{M}, t) + \rmi T_2(s + \tfrac{2\pi k}{M}, t) = \rme^{2\pi \rmi  k/M}(T_1(s, t) + \rmi X_2(T, t)),
    \cr
T_3(s + \tfrac{2\pi k}{M}, t) = T_3(s, t).
}
\end{equation}

\noindent Consequently, at a given $t$, $\X(s + 2\pi k/M, t)$ lay in the same orthogonal plane to the $z$-axis, for all $k\in\mathbb N$. This fact has important implications for implementing efficient numerical schemes, as we will see in Section \ref{s:numerical}.

Besides invariance by rotations, \eref{e:xt} is also mirror invariant: if $\X$ is a solution of \eref{e:xt}, so is $\tilde\X(s, t) = (-X_1(-s, t), X_2(-s, t), X_3(-s, t))^T$, etc. Therefore, if $\tilde\X(s, 0) = \X(s, 0)$, then that symmetry will be again preserved during the evolution, i.e, $\tilde\X(s, t) = \X(s, t)$, for all $t$, or, in other words, $\X(s, t)$ and $\X(-s, t)$ will be symmetric with respect to the plane containing the $z$-axis and the $y$-axis. In our case, a regular polygon with $M$ sides like \eref{e:xt0} has $2M$ of those mirror symmetries. An important corollary that will be used later is that $\X(2\pi/M, t) - \X(0, t)$ is a positive multiple of the vector $(1, 0, 0)^T$.

\subsection{Galilean invariance of the NLS equation}

\label{s:galilean}

One of the fundamental symmetries of the NLS equation \eref{e:schr} is the so-called Galilean transformations, i.e., if $\psi$ is a solution of \eref{e:schr}, so is
\begin{equation}
\label{e:galilean}
\tilde \psi_k(s, t) \equiv \rme^{\rmi ks - \rmi k^2t}\psi(s - 2kt, t), \quad\forall k, t\in\mathbb R.
\end{equation}

\noindent Therefore, if we choose an initial datum such that $\tilde \psi_k(s, 0) = \psi(s, 0)$, for all $k\in\mathbb R$, i.e., such that $\psi(s, 0) = \rme^{\rmi ks}\psi(s, 0)$, for all $k\in\mathbb R$, then, if the solution is unique, $\psi(s, t) = \rme^{\rmi ks - \rmi k^2t}\psi(s - 2kt, t)$, for all $k, t\in\mathbb R$. In this paper, this symmetry is enough for our purposes; for the full symmetry group of the NLS equation, see \cite[Section 17.2]{Ibragimov}.

In the case of a planar initial datum $\X(s, 0)$ for \eref{e:xt}, the torsion is equal to zero and, from \eref{e:ab}, $\psi(s, 0)$ is the curvature of $\X(s, 0)$. Hence, the $\psi(s, 0)$ corresponding to a regular polygon with $M$ sides is given by \eref{e:curvatureM}, i.e.,
\begin{equation}
\label{e:psis0}
\psi(s, 0) = \frac{2\pi}{M}\sum_{k = -\infty}^{\infty}\delta(s - \tfrac{2\pi k}{M}).
\end{equation}

\noindent Like \eref{e:curvatureM}, $\psi(s, 0)$ is $2\pi/M$-periodic and, since \eref{e:schr} is invariant with respect to space translations, $\psi(s, t)$ is also $2\pi/M$-periodic, for all $t\in\mathbb R$. Moreover, it satisfies $\psi(s, 0) = \rme^{\rmi Mks}\psi(s, 0)$, for all $k\in\mathbb Z$.
Therefore, the Galilean transformation holds, i.e., $\psi(s, t) = \rme^{\rmi Mks - \rmi(Mk)^2t}\psi(s - 2Mkt, t)$, for all $k\in\mathbb Z$, and for all $t\in\mathbb R$. Bearing in mind the previous arguments, we calculate the Fourier coefficients of $\psi(s, t)$:
\begin{eqnarray}
\label{e:Fourier}
\hat\psi(j, t) & = \frac{M}{2\pi}\int_0^{2\pi/M}\rme^{-\rmi Mjs}\psi(s, t)ds
    \cr
& = \frac{M}{2\pi}\int_0^{2\pi/M}\rme^{-\rmi Mjs}\left[\rme^{\rmi Mks - \rmi(Mk)^2t}\psi(s - 2Mkt, t)\right]ds
    \cr
& = \rme^{-\rmi (Mk)^2t - \rmi M(j-k)(2Mkt)}\hat\psi(j - k, t).
\end{eqnarray}

\noindent This identity holds for all $j$ and $k$. In particular, evaluating both sides at $j = k$, $\hat\psi(k, t) = \rme^{-\rmi (Mk)^2t}\hat\psi(0, t)$, so $\psi$ can be expressed as
\begin{equation}
\label{e:psist}
\psi(s, t) = \hat\psi(0, t)\sum_{k = -\infty}^\infty \rme^{-\rmi (Mk)^2t + \rmi Mks},
\end{equation}

\noindent where $\hat\psi(0, t)$ is a constant that depends on time and that has to be chosen in such a way that the corresponding curve $\X$ and tangent vector $\T$ are periodic of period $2\pi$ for a fixed $t$. Remark that, from the arguments following \eref{e:hasimoto2}, $\hat\psi(0, t)$ can be taken up to a constant with modulus one, so, in this paper, we assume without loss of generality that $\hat\psi(0, t)$ is real for any given $t$. When $t = 0$, we have trivially $\hat\psi(0, 0) = 1$; then, combining \eref{e:psis0} and \eref{e:psist}, we recover the following well-known identity:
\begin{equation}
\label{e:psiidentity}
\sum_{k = -\infty}^\infty \rme^{\rmi (Mk)s} \equiv \frac{2\pi}{M}\sum_{k = -\infty}^{\infty}\delta(s - \tfrac{2\pi k}{M}).
\end{equation}

\noindent Observe that, if we solve the free Schr\"odinger equation $\psi_t = \rmi \psi_{ss}$ for \eref{e:psis0}, we get
\begin{equation}
\label{e:psis1}
\psi(s, t) = \sum_{k = -\infty}^\infty \rme^{-\rmi (Mk)^2t + \rmi(Mk)s} = \theta\left(\frac{M}{2\pi}, \frac{M^2}{\pi}\right),
\end{equation}

\noindent where $\theta(s, t)$ is the well-known Jacobi theta function \eref{e:jacobi}. Equation \eref{e:psis1} is the same as \eref{e:psist}, but with $\hat\psi(0, t)\equiv 1$, for all $t\in\mathbb R$.

\subsection{$\psi(s, t)$ for rational multiples of $t = 2\pi/M^2$}

\label{s:rational}

Let us evaluate \eref{e:psist} at $t = t_{pq} = (2\pi/M^2)(p/q)$, where $p\in\mathbb Z$, $q\in\mathbb N$, and we can suppose without loss of generality that $\gcd(p, q) = 1$. Then, bearing in mind the identity \eref{e:psiidentity},
\begin{eqnarray}
\label{e:psistpq}
\psi(s, t_{pq}) & = \hat\psi(0, t_{pq})\sum_{k = -\infty}^\infty \rme^{-2\pi \rmi (p/q)k^2+\rmi Mks}
    \cr
& = \hat\psi(0, t_{pq})\sum_{m = 0}^{q - 1}\rme^{-2\pi \rmi (p/q)m^2 + \rmi Mms} \sum_{k = -\infty}^\infty \rme^{\rmi Mqks}
    \cr
& = \frac{2\pi}{Mq}\hat\psi(0, t_{pq})\sum_{m = 0}^{q - 1} \rme^{-2\pi \rmi (p/q)m^2 + \rmi Mms} \sum_{k = -\infty}^\infty \delta(s - \tfrac{2\pi k}{Mq})
    \cr
& = \frac{2\pi}{Mq}\hat\psi(0, t_{pq})\sum_{k = -\infty}^\infty\sum_{m = 0}^{q - 1}G(-p, m, q) \delta(s - \tfrac{2\pi k}{M} - \tfrac{2\pi m}{Mq}),
\end{eqnarray}

\noindent where $G(-p, m, q)$ was defined in \eref{e:Gabc}. Moreover, from \eref{e:absGabc}, it follows that
\begin{equation}
\label{e:G-pmq2}
G(-p, m, q) =
\cases{
\sqrt q \rme^{\rmi \theta_m}, & if $q$ odd,
    \\
\sqrt{2q} \rme^{\rmi \theta_m}, & if $q$ even $\wedge$ $q/2\equiv m\bmod 2$,
    \\
0, & if $q$ even $\wedge$ $q/2\not\equiv m\bmod 2$,
}
\end{equation}

\noindent for a certain angle $\theta_m$ that depends on $m$ (as well as, of course, $p$ and $q$). Therefore, $\psi(s, t)$ in $s\in[0,2\pi)$ has evolved from $M$ Dirac deltas at $t = 0$, to $Mq$ deltas at $t_{pq}$, for $q$ odd; and to $Mq/2$ deltas at $t_{pq}$, for $q$ even. Furthermore, introducing \eref{e:G-pmq2} into \eref{e:psistpq}, and restricting ourselves to $k = 0$, i.e., $s\in[0, 2\pi/M)$, we conclude that
\begin{equation}
\label{psistpq}
\psi(s, t_{pq})  =
\cases{
\frac{2\pi}{M\sqrt q}\hat\psi(0, t_{pq})\sum_{m = 0}^{q - 1}\rme^{\rmi \theta_m} \delta(s - \tfrac{2\pi m}{Mq}), & if $q$ odd,
    \\
\frac{2\pi}{M\sqrt {\tfrac{q}{2}}}\hat\psi(0, t_{pq})\sum_{m = 0}^{q/2 - 1}\rme^{\rmi \theta_{2m+1}} \delta(s - \tfrac{4\pi m + 2\pi}{Mq}), & if $q/2$ odd,
    \\
\frac{2\pi}{M\sqrt {\tfrac{q}{2}}}\hat\psi(0, t_{pq})\sum_{m = 0}^{q/2 - 1}\rme^{\rmi \theta_{2m}} \delta(s - \tfrac{4\pi m}{Mq}), & if $q/2$ even.
}
\end{equation}

\noindent The coefficients multiplying the Dirac deltas are in general not real, except for $t = 0$ and $t_{1,2} = \pi/M^2$. Therefore, $\psi(s, t_{pq})$ does not correspond to a planar polygon, but to a skew polygon with $Mq$ sides, for $q$ odd; and to a skew polygon with $Mq/2$ sides, for $q$ even. Moreover, $\hat\psi(0, t_{pq})$ has to be determined in such a way that the polygon is closed; this is done in the next section, where we show that its rather involved choice is unique, concluding that $\psi(s, t)$ is periodic in time, with period $2\pi / M^2$.

Since the Dirac deltas are equally spaced at $t = t_{pq}$, the length of the sides is the same. Observe that, if $q / 2$ is odd, i.e., if $q\equiv 2\bmod 4$, then there is no vertex at $s = 0$. For instance, if $p = 1$, $q = 2$, i.e., $t = \pi / M^2$, then
\begin{equation}
\psi(s, t_{1,2}) = \frac{2\pi}{M}\sum_{k = -\infty}^{\infty}\delta(s - \tfrac{2\pi (2k+1)}{2M}).
\end{equation}

\noindent In this case, the coefficients multiplying the Dirac deltas are real, so $\psi(s, t_{1, 2})$ is simply the curvature, and $\psi(0, t_{1, 2}) = 1$ has been chosen in order for the integral of $\psi(s, t_{1, 2})$ over $s\in[0, 2\pi)$ to be $2\pi$. Bearing in mind the symmetries of the problem, we conclude that the corresponding polygon is a regular planar polygon with $M$ sides, which lives at a certain plane $z = \mbox{const}$, but with the axes turned by $\pi / M$ degrees with respect to the initial datum, as R. Jerrard and D. Smets predicted in \cite{didier}.

\subsection{Recovering $\X$ and $\T$ from $\psi$ at $t = t_{pq}$}

\label{s:recovering}

Given a function $\psi(s, t_{pq}) = \alpha(s, t_{pq}) + \rmi\beta(s, t_{pq})$, recovering the basis vectors $\T$, $\pp$ and $\qq$ from $\psi$ implies integrating \eref{e:Te1e2ab}. As seen in the previous section, at $t = t_{pq}$, $\psi(s, t)$ is a sum of $Mq$ (if $q$ odd) or $Mq/2$ (if $q$ even) equally spaced Dirac deltas in $s\in[0,2\pi)$, that corresponds to a skew polygon $\X(s, t_{pq})$ of $Mq$ or $Mq/2$ sides. To integrate \eref{e:Te1e2ab}, we have to understand the transition from one side of the polygon to the next one. In order to do that, we reduce ourselves, without loss of generality, to a certain $\psi$ formed by one single Dirac delta located at $s = 0$, i.e, $\psi(s) = (a + \rmi b)\delta (s)$, so we have to integrate
\begin{equation}
\label{e:systemTe1e2}
\pmatrix{
\T
    \cr
\pp
    \cr
\qq
}_s
=
\pmatrix{
0 & a\delta(s) & b\delta(s)
    \cr
-a\delta(s) & 0 & 0
    \cr
-b\delta(s) & 0 & 0
}
\cdot
\pmatrix{
\T
    \cr
\pp
    \cr
\qq
} = \delta(s)\mathbf A\cdot\pmatrix{
\T
    \cr
\pp
    \cr
\qq
}.
\end{equation}

\noindent This system corresponds to some basis vectors $\T(s)$, $\pp(s)$ and $\qq(s)$ constant everywhere, except at $s = 0$, where there is a singularity. Let us suppose that, for $s < 0$, $\T(s) \equiv \T(0^-)$, $\pp(s) \equiv \pp(0^-)$, $\qq(s) \equiv \qq(0^-)$, so we want to calculate $\T(0^+)$, $\pp(0^+)$, $\qq(0^+)$, such that, for $s > 0$, $\T(s) \equiv \T(0^+)$, $\pp(s) \equiv \pp(0^+)$, $\qq(s) \equiv \qq(0^+)$. Then,
\begin{eqnarray}
\label{e:u+u-}
\left(
    \begin{array}{c}
\T(0^+)^T
    \cr
\hline
\pp(0^+)^T
    \cr
\hline
\qq(0^+)^T
    \end{array}
\right)
    & =
\exp\left[
\mathbf A
\int_{0^-}^{0^+}\delta(s')ds'\right]
    \cdot
\left(
    \begin{array}{c}
\T(0^-)^T
    \cr
\hline
\pp(0^-)^T
    \cr
\hline
\qq(0^-)^T
    \end{array}
\right)
    \cr
    & =
\exp(\mathbf A)\cdot
\left(
    \begin{array}{c}
\T(0^-)^T
    \cr
\hline
\pp(0^-)^T
    \cr
\hline
\qq(0^-)^T
    \end{array}
\right)
.
\end{eqnarray}

\noindent We can compute explicitly the matrix exponential. Writing $a + \rmi b\equiv \rho \rme^{\rmi \theta}$, then
\begin{eqnarray}
\label{e:expA}
\exp(\A) & =
\pmatrix{
c_\rho & s_\rho c_\theta & s_\rho s_\theta
    \cr
-s_\rho c_\theta & c_\rho c_\theta^2 + s_\theta^2 & [c_\rho - 1]c_\theta  s_\theta
    \cr
-s_\rho s_\theta & [c_\rho - 1]c_\theta s_\theta & c_\rho s_\theta^2 + c_\theta^2
},
\end{eqnarray}

\noindent where $c_\rho = \cos(\rho)$, $s_\rho = \sin(\rho)$, $c_\theta = \cos(\theta)$, $s_\theta = \sin(\theta)$. The matrix $\exp(\A)$ is a rotation matrix \cite{trenkler2008}; more precisely, it is the rotation matrix corresponding to a rotation about an axis $(0, \sin(\theta), -\cos(\theta))$ by an angle $\rho$.

Coming back to the general form of $\psi$, we have to integrate \eref{e:systemTe1e2} $Mq$ or $Mq/2$ times to obtain a closed skew, i.e., non-planar polygon with $Mq$ or $Mq/2$ sides. But, according to \eref{psistpq}, in $s\in[0, 2\pi/M)$,
\begin{equation}
\psi(s, t_{pq})  =
\cases{
\sum_{m = 0}^{q - 1}(\alpha_m + \rmi\beta_m) \delta(s - \tfrac{2\pi m}{Mq}), & if $q$ odd,
    \\
\sum_{m = 0}^{q/2 - 1}(\alpha_{2m+1} + \rmi\beta_{2m+1}) \delta(s - \tfrac{4\pi m + 2\pi}{Mq}), & if $q/2$ odd,
    \\
\sum_{m = 0}^{q/2 - 1}(\alpha_{2m} + \rmi\beta_{2m}) \delta(s - \tfrac{4\pi m}{Mq}), & if $q/2$ even;
}
\end{equation}

\noindent hence, we conclude that, at $t=t_{pq}$, the angle $\rho$ between two adjacent sides is constant.
\begin{equation}
\label{e:absolutevaluea_mb_m}
|\alpha_m + \rmi\beta_m| = \rho =
\cases{
\frac{2\pi}{M\sqrt q}\hat\psi(0, t_{pq}), & if $q$ odd,
    \cr
\frac{2\pi}{M\sqrt {\tfrac{q}{2}}}\hat\psi(0, t_{pq}), & if $q$ even $\wedge$ $q/2\equiv m\bmod 2$,
    \cr
0, & if $q$ even $\wedge$ $q/2\not\equiv m\bmod 2$.
}
\end{equation}

\noindent Furthermore, the structure of the polygon is completely determined by the quantities $\theta_m$ appearing in the generalized quadratic Gau{\ss} sum, where $\alpha_m + \rmi\beta_m = \rho \rme^{\rmi \theta_m}$. Let be $\M_m$ the rotation matrix corresponding to $(\alpha_m + \rmi\beta_m)\delta(s)$. If $\alpha_m + \rmi\beta_m \equiv 0$, $\M_m$ is simply an identity matrix and can be ignored. Otherwise, $\M_m$ is obtained after evaluating \eref{e:expA}, at $\theta = \theta_m$. Therefore,
\begin{eqnarray*}
\left(
    \begin{array}{c}
\T(\tfrac{2\pi(k+1)}{Mq}^-)^T
    \cr
\hline
\pp(\tfrac{2\pi(k+1)}{Mq}^-)^T
    \cr
\hline
\qq(\tfrac{2\pi(k+1)}{Mq}^-)^T
    \end{array}
\right)
& =
\left(
    \begin{array}{c}
\T(\tfrac{2\pi k}{Mq}^+)^T
    \cr
\hline
\pp(\tfrac{2\pi k}{Mq}^+)^T
    \cr
\hline
\qq(\tfrac{2\pi k}{Mq}^+)^T
    \end{array}
\right)
= \M_k\cdot
\left(
    \begin{array}{c}
\T(\tfrac{2\pi k}{Mq}^-)^T
    \cr
\hline
\pp(\tfrac{2\pi k}{Mq}^-)^T
    \cr
\hline
\qq(\tfrac{2\pi k}{Mq}^-)^T
    \end{array}
\right),
\end{eqnarray*}

\noindent or, equivalently,
\begin{equation}
\label{e:T2pik/mq}
\left(
    \begin{array}{c}
\T(\tfrac{2\pi k}{Mq}^+)^T
    \cr
\hline
\pp(\tfrac{2\pi k}{Mq}^+)^T
    \cr
\hline
\qq(\tfrac{2\pi k}{Mq}^+)^T
    \end{array}
\right)
    = \M_k\cdot\M_{k - 1}\cdot\ldots\M_1\cdot\M_0\cdot
\left(
    \begin{array}{c}
\T(0^-)^T
    \cr
\hline
\pp(0^-)^T
    \cr
\hline
\qq(0^-)^T
    \end{array}
\right),
\end{equation}

\noindent where $k\in\mathbb Z^+$, and $\M_k$ is periodic modulo $q$, i.e., $\M_{k + q} \equiv \M_k$. This formula is valid for both $q$ odd and $q$ even, bearing in mind that half of the $\M_k$ matrices are simply identity matrices if $q$ is even. This explains why the polygon has $q/2$ sides, when $q$ is even; equivalently, it could be regarded as having $q$ sides, but half of them being indistinguishable, because when the angle between two adjacent sides is zero, they merge into one single side.

In order for the polygon to be closed, we have to choose $\rho$ in \eref{e:absolutevaluea_mb_m} in such a way that the basis vectors $\T$, $\pp$ and $\qq$ are periodic, i.e.,
\begin{equation}
\left(
    \begin{array}{c}
\T(2\pi^-)^T
    \cr
\hline
\pp(2\pi^-)^T
    \cr
\hline
\qq(2\pi^-)^T
    \end{array}
\right)
=
\left(
    \begin{array}{c}
\T(0^-)^T
    \cr
\hline
\pp(0^-)^T
    \cr
\hline
\qq(0^-)^T
    \end{array}
\right),
\end{equation}

\noindent which is equivalent to imposing the condition
\begin{equation}
\label{e:productM}
\M_{Mq-1}\cdot \M_{Mq-2} \cdot \ldots \cdot \M_1 \cdot \M_0 \equiv \mathbf I.
\end{equation}

\noindent Let us define
\begin{equation}
\M = \M_{q - 1} \cdot \M_{q - 2} \cdot \ldots \cdot \M_{1} \cdot \M_0.
\end{equation}

\noindent From \eref{e:productM}, $\M$ is an $M$-th root of the identity matrix. Moreover, it is also a rotation matrix that induces a rotation of $2\pi/M$ degrees around a certain rotation axis. Therefore, we have to choose $\rho$ in order that any of the following properties is satisfied:
\begin{equation}
\Tr(\M) = 1 + 2\cos(\tfrac{2\pi}{M}), \qquad
\lambda(\M) = \{1, \rme^{2\pi \rmi  / M}, \rme^{-2\pi \rmi  / M}\},
\end{equation}

\noindent were $\Tr(\M)$ and $\lambda(\M)$ denote the trace and the spectrum of $\M$, respectively. We worked with the trace, because it is algebraically easier. In order to understand its structure, we analyzed a few cases ($q = 2, 3, 4, 5, 6, 8\ldots$) with the help of a symbolic manipulator. The obtained results strongly suggest that, for any $q$ (and for any $p$ coprime with it) the only possible choice for $\cos(\rho)$ is given by \eref{e:cosrho}.
Remark that if $q$ is odd, the angle $\rho$ between two adjacent sides is the same for $q$ and for $2q$.
Although finding a universal proof that \eref{e:cosrho} holds for any $q$ goes beyond the scope of this paper, we have checked it for a few more $q$. Moreover, this conjecture is in agreement with our numerical simulations, as we will see in Section \ref{s:experiments}. Therefore, we think that there is compelling evidence that \eref{e:cosrho} is valid for all $q\in\mathbb N$.

Once we have found the correct choice of $\rho$, we get from \eref{e:absolutevaluea_mb_m},
\begin{equation}
\label{e:hatpsi}
\hat\psi(0, t_{pq}) =
\cases{
\frac{M\sqrt q}{2\pi}\arccos\left(2\cos^{2/q}(\tfrac{\pi}{M}) - 1\right), & if $q$ odd,
    \\
\frac{M\sqrt {\tfrac{q}{2}}}{2\pi}\arccos\left(2\cos^{4/q}(\tfrac{\pi}{M}) - 1\right), & if $q$ even.
}
\end{equation}

\noindent It is straightforward to check that $\hat\psi(0, t_{1,2}) = 1$. Notice also that
\begin{equation}
\lim_{q\to\infty} \hat\psi(0, t_{pq}) = \frac{M\sqrt 2}{\pi}(-\ln(\cos(\tfrac{\pi}{M})))^{1/2},
    \qquad
\displaystyle\lim_{M\to\infty} \hat\psi(0, t_{pq}) = 1.
\end{equation}

\noindent Finally, from \eref{e:psistpq} and \eref{e:hatpsi}, we conclude that $\psi(s, t_{pq})$, $s\in[0,2\pi/M)$, is given by \eref{e:psistpq0}. This completes the proof of Theorem \ref{t:theorem}.

\subsection{Some extra information}

\label{s:rigid}

Knowing the function $\psi(s, t_{pq})$, the basis vectors $\T$ (and, hence, the curve $\X$), $\pp$ and $\qq$ can be completely determined up to a rigid movement, which, on the other hand, for an arbitrary polygon, can be pretty complicated to specify. However, since our initial datum is a regular planar polygon, the symmetries of the solution are very advantageous.

The easiest way to understand and to implement numerically the correct rotation (although not necessarily the simplest option for symbolic manipulation) is to work with $\X$. To that aim, given a rational time $t_{pq}$,
we compute first the associated rotation matrices $\M_m$. Then, by means of \eref{e:T2pik/mq}, we obtain up to a rotation the piecewise constant vectors $\T$, $\pp$ and $\qq$, which we denote $\tilde\T$, $\tpp$ and $\tqq$:
\begin{equation}
\left(
    \begin{array}{c}
\tilde\T(\tfrac{2\pi k}{Mq}^-)^T
    \cr
\hline
\tpp(\tfrac{2\pi k}{Mq}^-)^T
    \cr
\hline
\tqq(\tfrac{2\pi k}{Mq}^-)^T
    \end{array}
\right)
    =
\left(
    \begin{array}{c}
\tilde\T(\tfrac{2\pi k - 2\pi}{Mq}^+)^T
    \cr
\hline
\tpp(\tfrac{2\pi k - 2\pi}{Mq}^+)^T
    \cr
\hline
\tqq(\tfrac{2\pi k - 2\pi}{Mq}^+)^T
    \end{array}
\right),
\end{equation}

\noindent and the case $k = 0$ can be assigned any value, for instance the identity matrix. Again, it is not necessary to consider separately $q$ odd and $q$ even. Once $\tilde\T$ has been obtained, $\tilde\X$, which is $\X$ up to a rigid movement, is computed recursively:
\begin{equation}
\tilde\X(\tfrac{2\pi k + 2\pi}{Mq}) = \tilde\X(\tfrac{2\pi k}{Mq}) + \tfrac{2\pi}{Mq}\tilde\T(\tfrac{2\pi k}{Mq}^+),
\end{equation}

\noindent where, again, $\tilde\X(0)$ can be given any value, for instance $\tilde\X(0) = (0, 0, 0)^T$. To obtain the correct rotation for $\X$, we use the symmetries explained in Section \ref{s:formulation}. In particular, we use the fact that, at any time $t$, $\X(2\pi k/M)$, $k = 0, \ldots, M-1$, have to be coplanar and lay on a plane orthogonal to the $z$-axis; and that $\X(2\pi/M) - \X(0)$ is a positive multiple of $(1, 0, 0)^T$. An efficient algorithm is as follows:
\begin{enumerate}

\item Compute $\mathbf v^+ = \frac{\tilde\X(2\pi/M) - \tilde\X(0)}{\|\tilde\X(2\pi/M) - \tilde\X(0)\|}$ and $\mathbf v^- = \frac{\tilde\X(0) - \tilde\X(-2\pi/M)}{\|\tilde\X(0) - \tilde\X(-2\pi/M)\|}$.

\item Compute $\mathbf w = \mathbf v^- \wedge \mathbf v^+$.

\item Compute the scalar product between $\mathbf w$ and $(0, 0, 1)^T$, i.e., $w_3 = \mathbf w\cdot (0, 0, 1)^T$.

\item If $w_3 = 0$, $\mathbf R_1$ is the identity matrix. If not, $\mathbf R_1$ is the rotation matrix that induces a rotation of $\arccos(w_3)$ degrees around the axis given by the vector $\tfrac{\mathbf w\wedge(0, 0, 1)^T}{\|{\mathbf w\wedge(0, 0, 1)^T\|}}$.

\item Compute $\mathbf v^+_{new} = \mathbf R_1\cdot \mathbf v^+$.

\item Compute the rotation matrix $\mathbf R_2$ that induces a rotation of $\arccos(\mathbf v^+_{new}\cdot (1, 0, 0)^T)$ degrees around the $z$-axis.

\item Compute the sought rotation matrix, $\mathbf R = \mathbf R_2\cdot\mathbf R_1$.

\item Update $\T = \mathbf R\cdot\tilde\T$ and $\X = \mathbf R\cdot\tilde\X$.

\end{enumerate}

\noindent Once we have computed the correct rotation of $\X$ and $\T$, we still have to specify a translation for $\X$. Again, bearing in mind the symmetries of the problem, we translate $\X$ in such a way that its mass center, given by the mean of $\X(2\pi k/M)$, $k = 0, \ldots, M - 1$, is in the $z$-axis. Thus, the position of $\X$ is completely specified, up to a vertical translation.

\section{Numerical method}

\label{s:numerical}

In this section, we simulate numerically \eref{e:xt} and \eref{e:schmap} taking respectively \eref{e:xt0} and \eref{e:tt0} as the initial data. There have been a couple of papers devoted to reproducing numerically the self-similar solutions of the Schr\"odinger map equation \eref{e:schmap}. In \cite{buttke87}, a Crank-Nicholson scheme was consider, together with a finite-difference discretization of $\Tss$. Later on, in \cite{DelahozGarciaCerveraVega09}, both finite-difference discretizations and pseudo-spectral discretizations were considered to simulate \eref{e:schmap} and its equivalent on the hyperbolic plane. In particular, given the asymptotical structure of the self-similar solutions of \eref{e:schmap}, a truncated domain $s\in[-L, L]$, $L\gg1$, with a grid based on the Chebyshev nodes was found to be very convenient to represent those solutions. Nevertheless, due to the clustering of the Chebyshev nodes, an explicit scheme to advance in time implied the undesirable restriction $|\Delta t| = \mathcal O(1 / N^4)$, where $N$ is the number of nodes. In order to solve that, it was chosen to work with the stereographic projection of the tangent vector $\T=(T_1,T_2,T_3)$ over~$\mathbb C$,
\begin{equation}
\label{e:stereo}
z = x + \rmi y\equiv (x, y) \equiv
\left({\frac{T_1}{1 + T_3}}, {\frac{T_2}{1 + T_3}}\right),
\end{equation}

\noindent which transforms \eref{e:schmap} into a nonlinear Schr\"odinger equation:
\begin{equation}
\label{e:z_t} z_t = \rmi z_{ss} - \frac{2\rmi\bar z}{1 + |z|^2}z_s^2.
\end{equation}

\noindent The advantage of this equation, as opposed to \eref{e:schmap}, is that the higher-order term $z_{ss}$ can be treated implicitly, therefore eliminating or at least reducing significantly the restrictions on $\Delta t$. Indeed, working with \eref{e:z_t} was a very adequate choice for the purposes of \cite{DelahozGarciaCerveraVega09}, where a second-order semi-implicit backward differentiation formula was chosen. This scheme is very stable because it imposes a very strong decay in the high frequency modes; moreover, its low order is easily and effectively compensated by the use of an adaptive method, both in space and in time.

However, in this paper, working with \eref{e:z_t} does not seem to be such a good option. Indeed, unlike in \cite{DelahozGarciaCerveraVega09}, where an extremely high accuracy was required only for $0 < t\ll 1$, we are now interested in the behaviour of $\T$ at all times, for which a second-order scheme seems a very poor option. After many unsuccessful attempts, we admit that we have been unable to find a good higher-order scheme for \eref{e:z_t}. Furthermore, the stereographic projection \eref{e:stereo} may introduce an artificial singularity at the points near the south pole of the sphere. Finally, we are also interested in the evolution of the curve $\X$, for which working with \eref{e:z_t} is of limited help. Bearing in mind the previous arguments, we have opted to work directly with \eref{e:schmap} or, more precisely, with a combination of \eref{e:xt} and \eref{e:schmap}:
\begin{equation}
\label{e:XTt}
\cases{
\Xt = \T\wedge\Ts,
    \\
\Tt = \T\wedge\Tss,
}
\end{equation}

\noindent to which we have applied a fourth-order Runge-Kutta scheme in time:
\begin{eqnarray}
\label{e:rungekutta}
\mathbf A_X & = \T^{(n)} \wedge \Ts^{(n)}, \quad \mathbf A_T = \T^{(n)} \wedge \Tss^{(n)}, \quad
\T^{(A)} = \T^{(n)} + \tfrac{\Delta t}{2} \mathbf A_T,
    \cr
\mathbf B_X & = \T^{(A)} \wedge \Ts^{(A)}, \quad \mathbf B_T = \T^{(A)} \wedge \Tss^{(A)}, \quad
\T^{(B)} = \T^{(n)} + \tfrac{\Delta t}{2} \mathbf B_T,
    \cr
\mathbf C_X & = \T^{(B)} \wedge \Ts^{(B)}, \quad \mathbf C_T = \T^{(B)} \wedge \Tss^{(B)}, \quad
\T^{(C)} = \T^{(n)} + \Delta t \mathbf C_T,
    \cr
\mathbf D_X & = \T^{(C)} \wedge \Ts^{(C)}, \quad \mathbf D_T = \T^{(C)} \wedge \Tss^{(C)},
    \cr
\X^{(n + 1)} & = \X^{(n)} + \tfrac{\Delta t}{6}(\mathbf A_X + 2\mathbf B_X + 2\mathbf C_X + \mathbf D_X),
    \cr
\tilde\T & = \T^{(n)} + \tfrac{\Delta t}{6}(\mathbf A_T + 2\mathbf B_T + 2\mathbf C_T + \mathbf D_T), \quad \T^{(n + 1)} = \frac{\tilde\T}{\|\tilde\T\|}.
\end{eqnarray}

\noindent The approximations $\X^{(n)}$ and $\T^{(n)}$ refer to the numerically obtained values of $\X$ and $\T$ at $t^{(n)} \equiv n\Delta t$; the last line guarantees that $\T^{(n)}\in\mathbb S^2$, for all $n$. Additionally, we can also project $\T^{(A)}$, $\T^{(B)}$, $\T^{(C)}$ onto the unit sphere, but without significant improvement in the numerical results.

\subsection{Spatial symmetries and DFT}

\label{s:symmetries}

We have combined \eref{e:rungekutta} with a pseudo-spectral discretization directly in space. More precisely, since we deal with periodic solutions in $s\in[0,2\pi)$, we have simulated the evolution of the curve $\X$ and of the tangent vector $\T$ at $N$ equally spaced nodes $s_j = 2\pi j / N$, $j = 0, \ldots, N-1$. In order to compute $\Ts$, $\Tss$, etc., we remember that, given a periodical function $f(s)$ evaluated at $s_j$, its derivatives at $s_j$ can be spectrally approximated as
\begin{equation}
\label{e:fsfss}
f_s(s_j) = \sum_{k = -N/2}^{N/2-1}\rmi k\hat f(k)\rme^{2\pi \rmi jk/N}, \quad
f_{ss}(s_j) = -\sum_{k = -N/2}^{N/2-1}k^2\hat f(k)\rme^{2\pi \rmi jk/N},
\end{equation}

\noindent where
\begin{equation}
\label{e:hatfk}
\hat f(k) = \sum_{j = 0}^{N-1}f(s_j)\rme^{-2\pi \rmi jk/N}.
\end{equation}

\noindent Both \eref{e:fsfss} and \eref{e:hatfk}, which are respectively inverse and direct discrete Fourier transforms (DFT) of $N$ elements, can be computed efficiently by means of the Fast Fourier Transform (FFT) algorithm \cite{FFTW}. Moreover, since $\X$ and $\T$ are invariant with respect to rotations of $2\pi/M$ degrees around the $z$-axis, as shown in \eref{e:rotationXT}, we can reduce the computational cost of \eref{e:fsfss} and \eref{e:hatfk} to a DFT of $N/M$ elements. In what follows, we explain how to apply this idea to $\T$. This is valid, with no change, for any vector satisfying \eref{e:rotationXT}, i.e., $\X$, $\Ts$, $\Tss$, etc. Now, let us consider the first two components of $\T$. Denoting $Z(s_j) = T_1(s_j, t) + \rmi T_2(s_j, t)$, and bearing in mind that $Z(s_{j + N/M}) = Z(s_j)\rme^{2\pi\rmi/M}$,
\begin{eqnarray}
\hat Z(k) & = \sum_{j = 0}^{N - 1}Z(s_j)\rme^{-2\pi \rmi jk/N}
    \cr
& = \left(\sum_{l = 0}^{M - 1}\rme^{-2\pi \rmi l(k-1)/M}\right)\sum_{j = 0}^{N/M - 1}Z(s_j)\rme^{-2\pi \rmi jk/N}
    \cr
& =
\cases{
M\displaystyle{\sum_{j = 0}^{N/M - 1}}Z(s_j)\rme^{-2\pi \rmi jk/N}, & if $k \equiv 1 \bmod M$,
    \\
0, & if $k \not\equiv 1 \bmod M$.
}
\end{eqnarray}

\noindent Therefore, the non-zero $\hat Z(k)$ are obtained by a DFT of $N/M$ elements:
\begin{equation}
\label{e:ZMk}
\hat Z(Mk + 1) = M\sum_{j = 0}^{N/M - 1}\left[\rme^{-2\pi \rmi j/N}Z(s_j)\right]\rme^{-2\pi \rmi jk/(N/M)}.
\end{equation}

\noindent For the third component of $\T$, bearing in mind that $T_3(s_{j + N/M}) = T_3(s_j)$, we conclude
\begin{eqnarray}
\hat T_3(k) =
\cases{
M\displaystyle{\sum_{j = 0}^{N/M - 1}}T_3(s_j)\rme^{-2\pi \rmi jk/N}, & if $k \equiv 0 \bmod M$,
    \\
0, & if $k \not\equiv 0 \bmod M$,
}
\end{eqnarray}

\noindent so the non-zero $\hat T_3(k)$ are obtained again by a DFT of $N/M$ elements:
\begin{equation}
\label{e:T3Mk}
\hat T_3(Mk) = M\sum_{j = 0}^{N/M - 1}T_3(s_j)\rme^{-2\pi \rmi jk/(N/M)};
\end{equation}

\noindent since $T_3(s)$ is real, \eref{e:T3Mk} can be further simplified to a DFT of $N/(2M)$ elements.

\subsection{Stability and computational cost}

\label{s:cost}

Combining \eref{e:rungekutta} with a pseudo-spectral discretization in space, we obtain a scheme with a time-step restriction on $\Delta t$ that appears to be of the form $\Delta t \le C/N^2$, with $C \approx 11.3$. However, for a fixed $M$, we only want to simulate the evolution of the curve $\X$ and of the tangent vector $\T$ for $t\in[0, 2\pi/M^2]$. Denoting $\Delta t_{max}$ the biggest $\Delta t$ that makes \eref{e:rungekutta} stable, the smallest number of time-steps $N_t$ needed to reach $t = 2\pi / M^2$ satisfies
\begin{equation}
\label{e:timesteps}
N_t \ge \left\lceil\frac{2\pi/M^2}{\Delta t_{max}}\right\rceil = \left\lceil\frac{2\pi/M^2}{C/N^2}\right\rceil = \left(\frac{N}{M}\right)^2\left\lceil\frac{2\pi}{C}\right\rceil.
\end{equation}

\noindent From \eref{e:timesteps}, together with \eref{e:ZMk} and \eref{e:T3Mk}, it follows that, if the number of spatial nodes $N$ is of the form $N = 2^rM$, with $r\in\mathbb N$, then, the computational cost for simulating \eref{e:rungekutta} is exactly the same for any number of initial sides $M$. This important fact enables us to make consistent comparisons between different $M$, because we are considering an equivalent spatial resolution, i.e, we are using the same number of points $N/M$ per initial side. We take $N/M$ a power of two to take the greatest possible advantage of the FFT. Obviously, for $N/M$ constant, $\Delta t$ decreases as $\mathcal O(1/M^2)$; therefore, the most accurate results will be expected for the largest $M$, as we will observe in our numerical experiments.

In the following section, we simulate $\X$ and $\T$ for initial data with different $M$. As we will see, despite its simplicity, \eref{e:rungekutta} gives surprisingly good results.

\section{Numerical experiments}

\label{s:experiments}

In Section \ref{s:Xt}, in order to construct algebraically the evolution of a regular polygon at rational times, we have done some very strong assumptions, the most important one being uniqueness. However, as we will see in this section, the numerical experiments are in complete agreement with their theoretical predictions. To illustrate this, we have simulated the evolution of the curve $\X$ and of the tangent vector $\T$ by means of \eref{e:rungekutta}, taking regular polygons with different numbers of sides $M$ as the initial data, and making $\X$ and $\T$ evolve until $t = 2\pi / M^2$. The initial data $\X(s, 0)$ and $\T(s, 0)$ are given respectively by \eref{e:xt0} and \eref{e:tt0}; in the case of $\X$, the non-vertex points are computed by linear interpolation. We have taken $N$ equally spaced nodes $s_j = 2\pi j / N\in[0, 2\pi)$, $j = 0, \ldots, N-1$; nevertheless, from the symmetries of $\X$ and $\T$, as explained in the previous section, we only have to describe the evolution of $\X$ and $\T$ at the first $N/M$ nodes $s_j\in[0, 2\pi/M)$, $j = 0, \ldots, N/M-1$. We have divided the time-interval $[0, 2\pi/M^2]$ in $N_t$ equally spaced time steps of length $\Delta t = (2\pi / M^2) / N_t$. Therefore, during the simulation of \eref{e:rungekutta}, we obtain $\X^{(n)}$ and $\T^{(n)}$, at $t^{(n)}$, $n = 1, \ldots, N_t$. In our experiments, bearing in mind \eref{e:timesteps}, we have chosen $N_t = 151200\cdot4^r$, for $N/M = 512\cdot2^r$.

Remember that, in Section \ref{s:Xt}, we had constructed algebraically the curve $\X$ up to a vertical movement. Therefore, in order to completely specify $\X$ at a given time, we would need to give the height of one point or, more conveniently, the height $h(t)$ of the mass center, which is precisely the mean of all the values $X_3(s_j, t)$:
\begin{equation}
h(t) = \frac{1}{N}\sum_{j = 0}^{N-1}X_3(s_j, t) = \frac{M}{N}\sum_{j = 0}^{N/M-1}X_3(s_j, t).
\end{equation}

\noindent We have observed that $h(t)$ can be approximated with very high accuracy by means of a constant multiplied by $t$; more precisely,
\begin{equation}
h(t) \approx \frac{h(2\pi / M^2)}{2\pi / M^2}t = c_Mt,
\end{equation}

\noindent where $c_M \equiv h(2\pi/M^2) / (2\pi/M^2)$ is the mean speed. In Table \ref{t:speed}, we give the maximum discrepancy between $h(t)$ and its linear approximation $c_Mt$, i.e., $\max_n|h(t^{(n)}) - c_Mt^{(n)}|$. We have considered different numbers of initial sides $M$ and different numbers of nodes $N$; the errors are very small and they seem to decrease as $\mathcal O(1/N)$. This gives very strong evidence that $h(t)$ is linear of the form $h(t) = c_Mt$. Observe that $c_M$, which is also offered in Table \ref{t:speed} for $N/M = 8192$, grows with $M$, tending to $1$, as $M$ tends to infinity, i.e. as the initial polygon $\X(s, 0)$ tends to a circle.
\begin{table}[t!]
\centering
\begin{tabular}{|c|c|c|c|c|c|c|}
\hline $M$ & $N/M=512$ & $N/M = 1024$ & $N/M = 2048$ & $N/M = 4096$ & $N/M = 8192$ & $c_M$
\\
\hline $3$ & $4.3096\cdot10^{-5}$ & $2.1206\cdot10^{-5}$ & $1.0886\cdot10^{-5}$ & $5.7953\cdot10^{-6}$ & $3.1123\cdot10^{-6}$ & $0.7644$
\\
\hline $4$ & $1.2398\cdot10^{-6}$ & $6.1344\cdot10^{-6}$ & $3.2140\cdot10^{-6}$ & $1.7316\cdot10^{-6}$ & $9.4280\cdot10^{-7}$ & $0.8826$
\\
\hline $5$ & $4.8504\cdot10^{-6}$ & $2.4191\cdot10^{-6}$ & $1.2807\cdot10^{-6}$ & $6.9338\cdot10^{-7}$ & $3.7928\cdot10^{-7}$ & $0.9286$
\\
\hline $6$ & $2.2848\cdot10^{-6}$ & $1.1441\cdot10^{-6}$ & $6.0905\cdot10^{-7}$ & $3.3044\cdot10^{-7}$ & $1.8113\cdot10^{-7}$ & $0.9517$
\\
\hline $7$ & $1.2167\cdot10^{-6}$ & $6.1060\cdot10^{-7}$ & $3.2607\cdot10^{-7}$ & $1.7710\cdot10^{-7}$ & $9.7195\cdot10^{-8}$ & $0.9650$
\\
\hline $8$ & $7.0721\cdot10^{-7}$ & $3.5594\cdot10^{-7}$ & $1.9014\cdot10^{-7}$ & $1.0333\cdot10^{-7}$ & $5.6754\cdot10^{-8}$ & $0.9735$
\\
\hline $9$ & $4.3905\cdot10^{-7}$ & $2.2140\cdot10^{-7}$ & $1.1828\cdot10^{-7}$ & $6.4303\cdot10^{-8}$ & $3.5336\cdot10^{-8}$ & $0.9792$
\\
\hline $10$ & $2.8697\cdot10^{-7}$ & $1.4489\cdot10^{-7}$ & $7.7407\cdot10^{-8}$ & $4.2093\cdot10^{-8}$ & $2.3139\cdot10^{-8}$ & $0.9832$
\\
\hline
\end{tabular}
\caption{$\max_n|h(t^{(n)}) - c_Mt^{(n)}|$, for different numbers of initial sides $M$, and of nodes $N$. $t^{(n)} = n\Delta t$; $c_M \equiv h(2\pi/M^2) / (2\pi/M^2)$. The values corresponding to $c_M$ have been calculated with $N/M = 8192$.} \label{t:speed}
\end{table}

In order to compare the values of $\X$ obtained numerically, which we label $\X_{num}$, with those obtained algebraically, which we label $\X_{alg}$, we have subtracted the vertical position of the center of mass from $\X_{num}$, i.e., we have analyzed the agreement between $\X_{num} - c_Mt(0,0,1)^T$ and $\X_{alg}$. More precisely, we have computed $\max_m(\max_j\|\X_{num}(s_j, t^{(m)}) - c_Mt^{(m)}(0, 0, 1)^T - \X_{alg}(s_j, t^{(m)})\|)$, where $s_j = 2\pi j/N$, $j = 0,\ldots, N/M-1$; $\|\cdot\|$ denotes the Euclidean distance; and $t^{(m)} = (2\pi/M^2)(m/1260)$, $m = 0, \ldots, 1260$. Notice that, in Section \ref{s:Xt}, we have constructed $\X_{alg}(s, t)$ only at the vertices, so the non-vertex points $\X_{alg}(s_j, t)$ are computed by linear interpolation. Observe also that comparing $\X_{num}$ and $\X_{alg}$ at $N_t+1$ time-instants would have been computationally unrealistic. Instead, we have chosen $1260 + 1$ time-instants, because $1260$ is a number large enough for our purposes, and with a convenient factorization, $1260 = 2^2\cdot3^2\cdot5\cdot7$.

The maximum value of $\|\X_{num}(s_j, t^{(m)}) - c_Mt^{(m)}(0, 0, 1)^T - \X_{alg}(s_j, t^{(m)})\|$ is given in Table \ref{t:errors}, for different $N/M$ and $M$. Bearing in mind that we are comparing $\X_{num}$ and $\X_{alg}$ globally for a large number of nodes and time-instants, and that $\max\|\X_{alg}\| > 1$, for all $M$, the results are, in our opinion, very remarkable, and strongly suggest that there is convergence between both approaches, as $M$ tends to infinity; moreover, since $\X_{alg}$ is periodical in time, with time-period $2\pi / M^2$, they also suggest that $\X$ is periodical in time with that period, up to a vertical movement with constant velocity. All this is also supported by the fact that the errors decrease as $M$ increases. Indeed, from \eref{e:timesteps}, we are taking the same number of time-steps $N_t$ for each $M$, but $\Delta t = (2\pi/M^2) / N_t$, i.e., $\Delta t = \mathcal O(1/M^2)$, which explains those smaller values for the last rows of Table \ref{t:errors} and also of Table \ref{t:speed}.

\begin{table}[t!]
\centering
\begin{tabular}{|c|c|c|c|c|c|}
\hline $M$ & $N/M=512$ & $N/M = 1024$ & $N/M = 2048$ & $N/M = 4096$ & $N/M = 8192$
\\
\hline $3$ & $2.4847\cdot10^{-3}$ & $1.3841\cdot10^{-3}$ & $8.1211\cdot10^{-4}$ & $4.9718\cdot10^{-4}$ & $3.0091\cdot10^{-4}$
\\
\hline $4$ & $1.1221\cdot10^{-3}$ & $6.9665\cdot10^{-4}$ & $4.2717\cdot10^{-4}$ & $2.5917\cdot10^{-4}$ & $1.7505\cdot10^{-4}$
\\
\hline $5$ & $6.8414\cdot10^{-4}$ & $4.2545\cdot10^{-4}$ & $2.6125\cdot10^{-4}$ & $1.5874\cdot10^{-4}$ & $1.1378\cdot10^{-4}$
\\
\hline $6$ & $4.6057\cdot10^{-4}$ & $2.8717\cdot10^{-4}$ & $1.7670\cdot10^{-4}$ & $1.0754\cdot10^{-4}$ & $7.9642\cdot10^{-5}$
\\
\hline $7$ & $3.3170\cdot10^{-4}$ & $2.0724\cdot10^{-4}$ & $1.2772\cdot10^{-4}$ & $7.7832\cdot10^{-5}$ & $5.8787\cdot10^{-5}$
\\
\hline $8$ & $2.5059\cdot10^{-4}$ & $1.5680\cdot10^{-4}$ & $9.6744\cdot10^{-5}$ & $5.9010\cdot10^{-5}$ & $4.5144\cdot10^{-5}$
\\
\hline $9$ & $1.9616\cdot10^{-4}$ & $1.2288\cdot10^{-4}$ & $7.5878\cdot10^{-5}$ & $4.6313\cdot10^{-5}$ & $3.5743\cdot10^{-5}$
\\
\hline $10$ & $1.5782\cdot10^{-4}$ & $9.8943\cdot10^{-5}$ & $6.1137\cdot10^{-5}$ & $3.7334\cdot10^{-5}$ & $2.8994\cdot10^{-5}$
\\
\hline
\end{tabular}
\caption{$\max_m(\max_j\|\X_{num}(s_j, t^{(m)}) - c_Mt^{(m)}(0, 0, 1)^T - \X_{alg}(s_j, t^{(m)})\|)$, where $s_j = 2\pi j / N$, $j = 0,\ldots, N/M-1$; $\|\cdot\|$ denotes the Euclidean distance; and $t^{(m)} = (2\pi / M^2)(m / 1260)$, $m = 0, \ldots, 1260$.} \label{t:errors}
\end{table}

\begin{figure}[t!]
\center
\includegraphics[width=0.5\textwidth, clip=true]{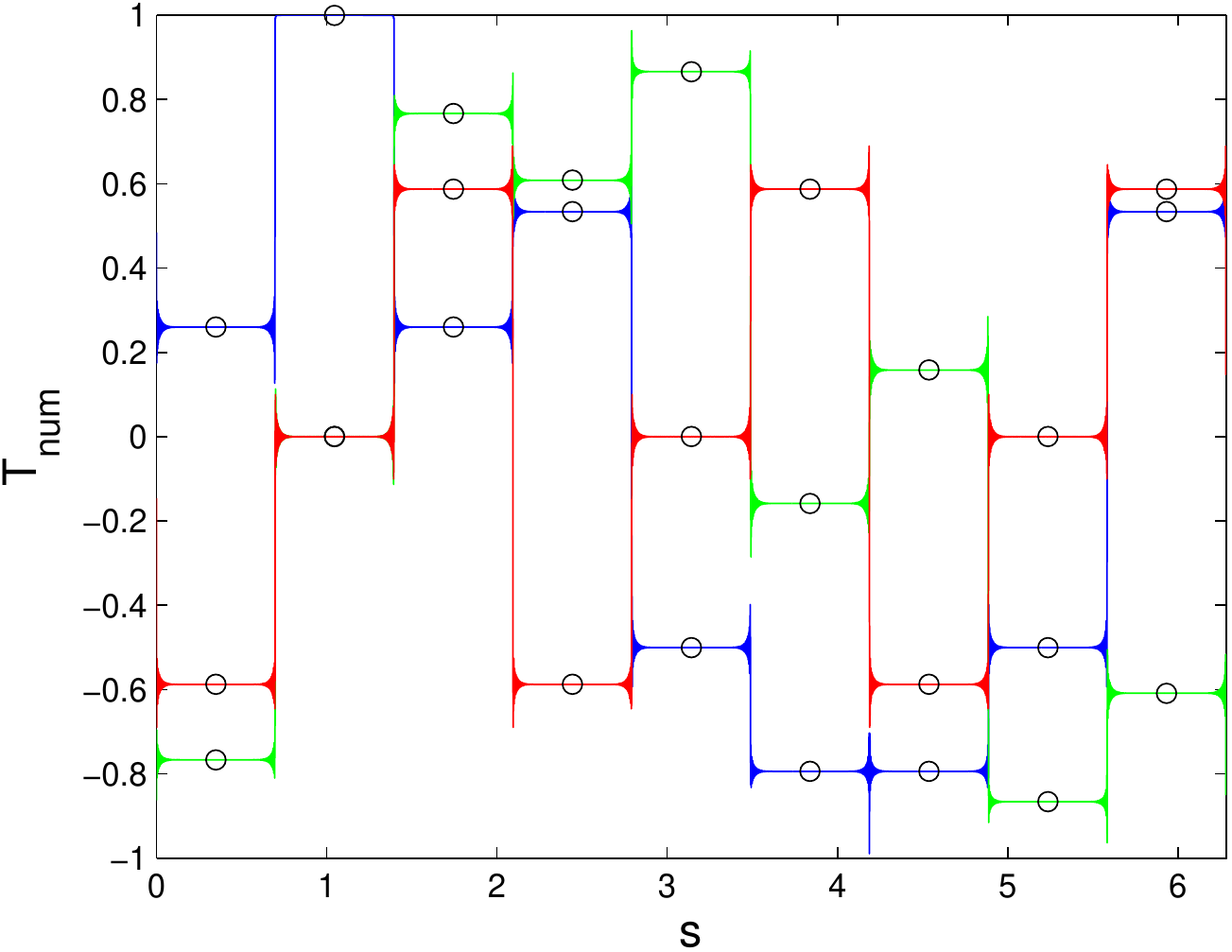}~\includegraphics[width=0.5\textwidth, clip=true]{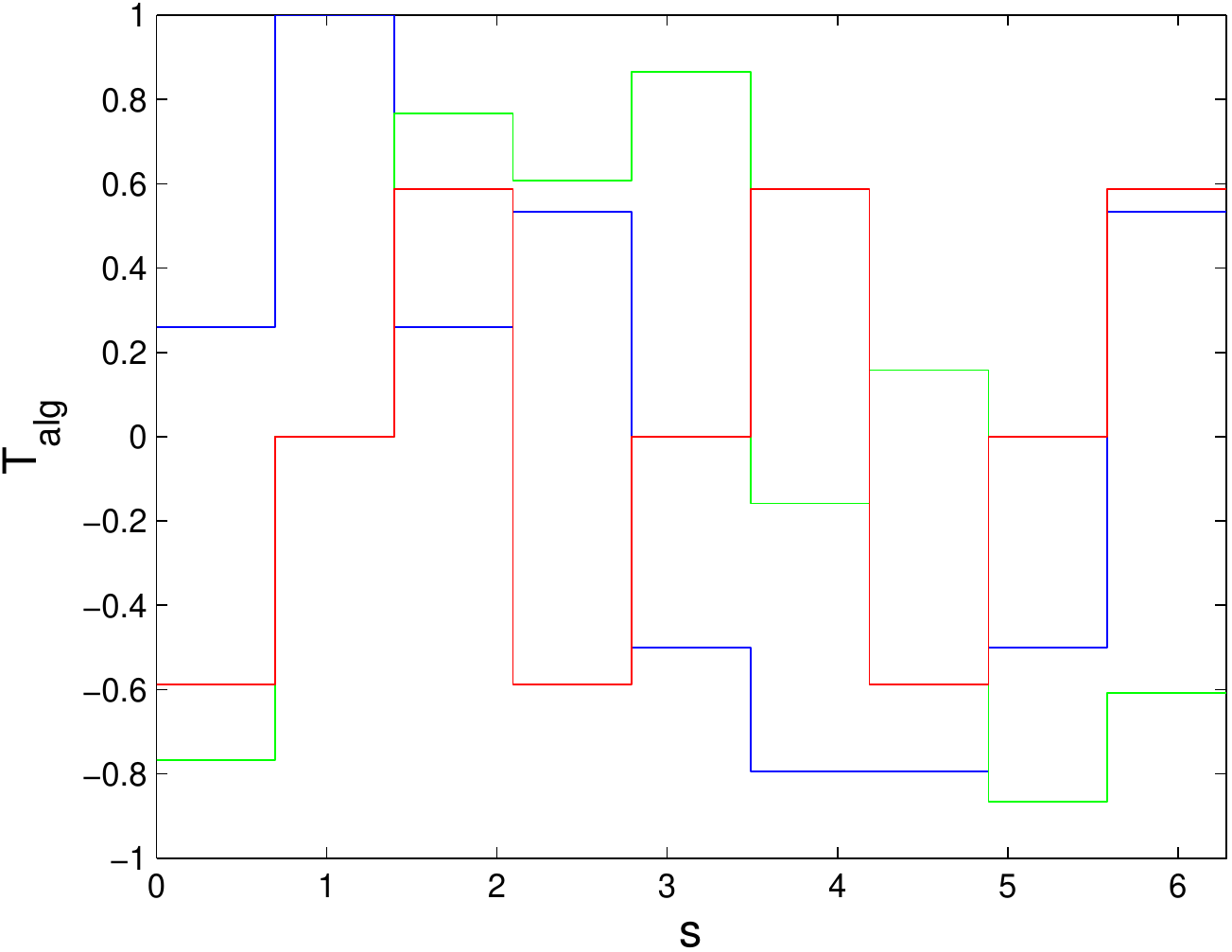}
\caption{$\T_{num}$ versus $\T_{alg}$, for $M = 3$, at $t_{1, 3} = 2\pi/27$. $T_1$ appears in blue, $T_2$ in green, $T_3$ in red. In $\T_{num}$, the Gibbs phenomenon is clearly visible. The black circles denote the points chosen for the comparisons.}\label{f:TM3p1q3}
\end{figure}

On the other hand, comparing the algebraically constructed $\T_{alg}$ with the numerically obtained $\T_{num}$ is more problematic. For example, in Figure \ref{f:TM3p1q3}, we have compared $\T_{num}$ (left) with $\T_{alg}$ (right), for an initial triangle, $M = 3$, at $t_{1, 3} = 2\pi/27$, $N/M = 8192$. The exact value of $\T$ at that time is given by $\T_{alg}$, which, by construction, is a piecewise constant function with exactly $Mq = 9$ pieces, that we denote $\T_i$, $i = 1, \ldots, 9$, and whose explicit values, in this case, have an easy algebraic expression:
\begin{equation}
\label{e:Texact}
\T_1 =
\pmatrix{
\sqrt[3]2-1 \cr -\sqrt{\sqrt[3]4-1} \cr 1 - \sqrt[3]4
},
\quad
\T_2 =
\pmatrix{
1 \cr 0 \cr 0
},
\quad
\T_3 =
\pmatrix{
\sqrt[3]2-1 \cr \sqrt{\sqrt[3]4-1} \cr \sqrt[3]4 - 1
};
\end{equation}

\noindent $\T_4$, $\T_5$ and $\T_6$ are respectively $\T_1$, $\T_2$ and $\T_3$ rotated $2\pi/3$ degrees around the $z$-axis; $\T_7$, $\T_8$ and $\T_9$ are respectively $\T_1$, $\T_2$ and $\T_3$ rotated $4\pi/3$ degrees around the $z$-axis. Observe that $\T_{num}$ is almost identical to $\T_{alg}$, but clearly exhibits a Gibbs-type phenomenon. However, in $\T_{num}$, if we take the central values of each interval (altogether 27 circles, indicated in Figure \ref{f:TM3p1q3} with a black circle), and compare them componentwise with \eref{e:Texact}, we obtain a maximum error equal to $3.3976\cdot10^{-10}$.

\subsection{$\X(0, t)$ and Riemann's non-differentiable function}

\label{s:riemann}

From the mirror symmetries of the initial data, we conclude that the curve $\X(\pi k/M, t)$, $k = 0, \ldots, 2M-1$, lives in a plane that contains the $z$-axis. Indeed, $s = \pi k/M$ are the only values for which $\X(s, t)$ describe a planar curve. In what follows, we will describe $\X(0, t)$, although all said here is immediately applicable to any $s = \pi k/M$. Since $\X(0, t) = (X_1(0, t), X_2(0, t), X_3(0, t))$ is planar, bearing in mind that $X_1(0, t) < 0$, $X_2(0, t) < 0$, we rotate $\X(0, t)$ clockwise $\pi/2-\pi/M$ degrees around the $z$-axis, until it lays on the plane $OYZ$, which we identify with $\mathbb C$. Then $\X(0, t)$ becomes
\begin{equation}
z(t) = -\|(X_1(0, t), X_2(0, t))\| + \rmi X_3(0, t).
\end{equation}

\begin{figure}[t!]
\center
\includegraphics[width=0.3333\textwidth, clip=true]{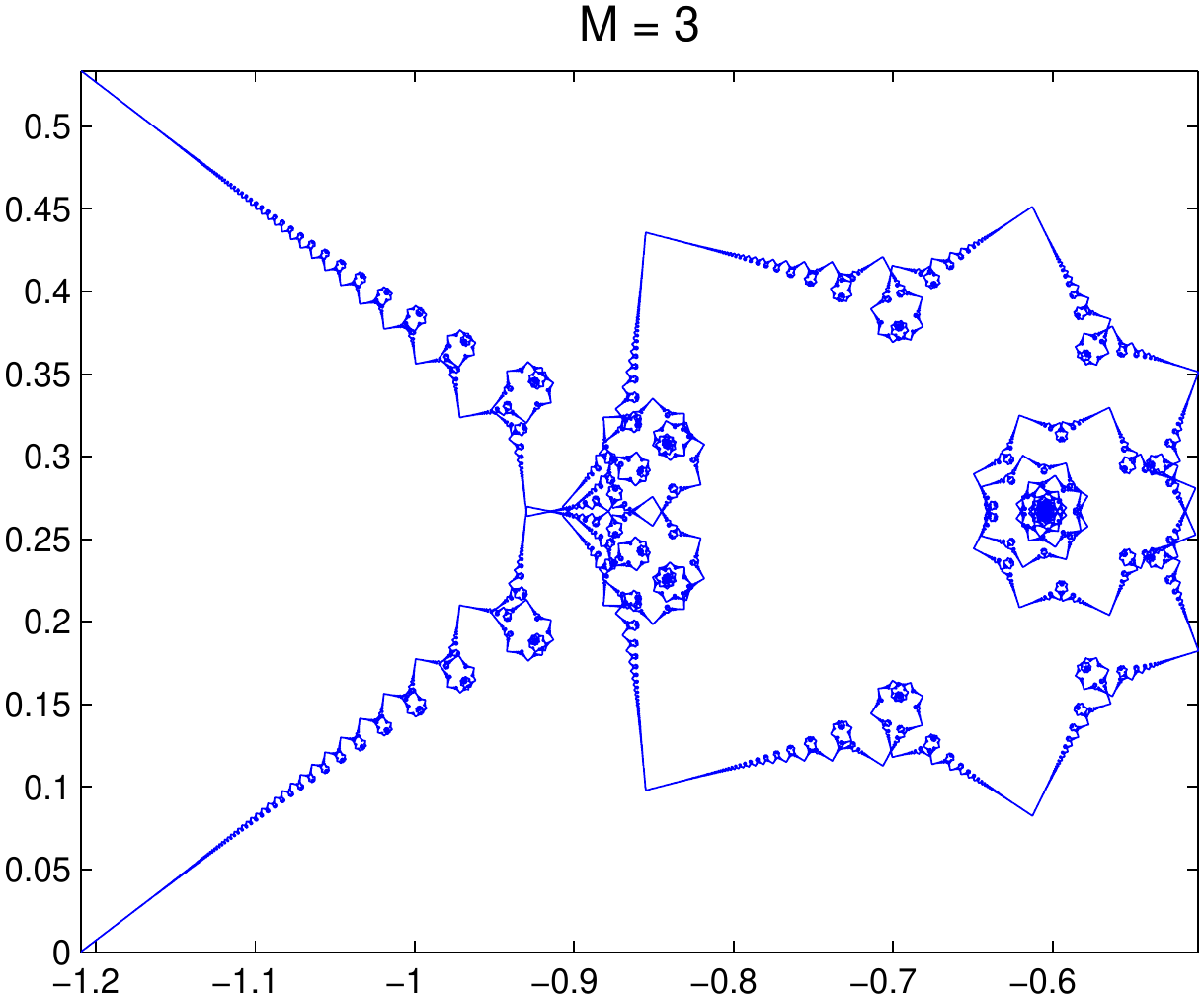}~
\includegraphics[width=0.3333\textwidth, clip=true]{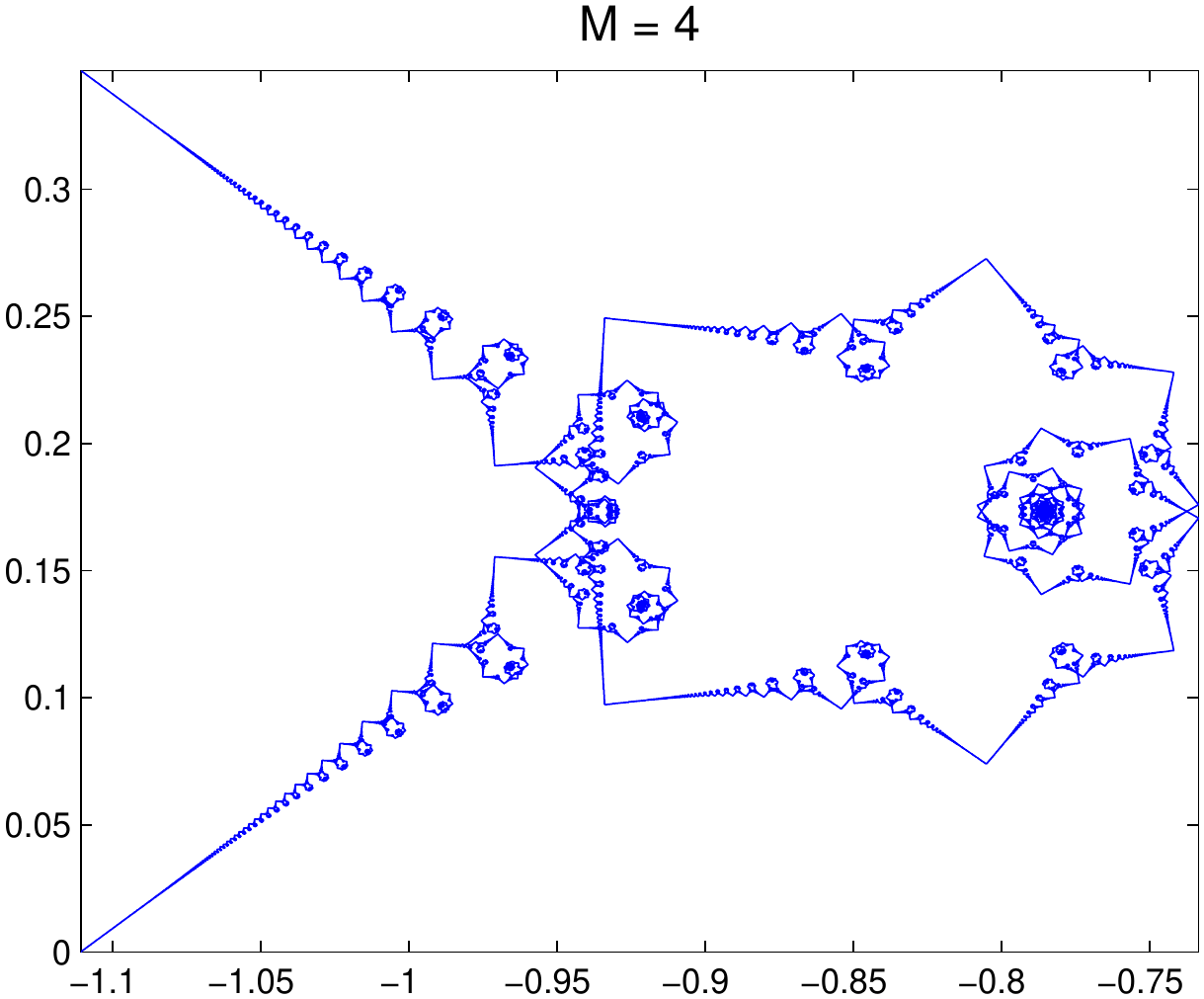}~
\includegraphics[width=0.3333\textwidth, clip=true]{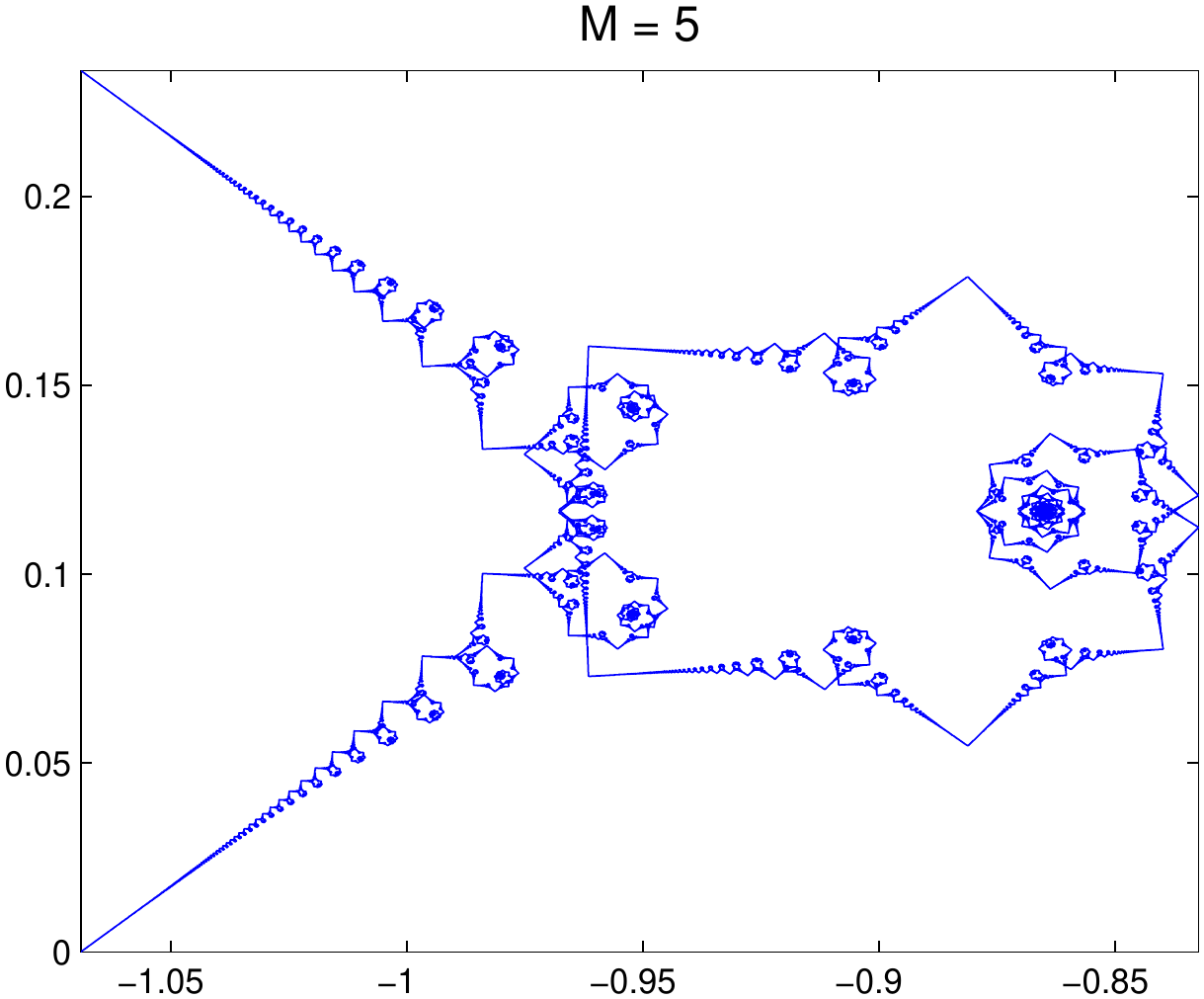}

\caption{$z(t)$, for $M = 3, 4, 5$. $N/M = 8192$. We have plotted the $z(t)$ corresponding to all the $N_t + 1$ points $\X^{(n)}(0)$ obtained during the execution of \eref{e:rungekutta}.}\label{f:Xs0}
\end{figure}

\begin{figure}[t!]
\center

\includegraphics[width=0.3333\textwidth, clip=true]{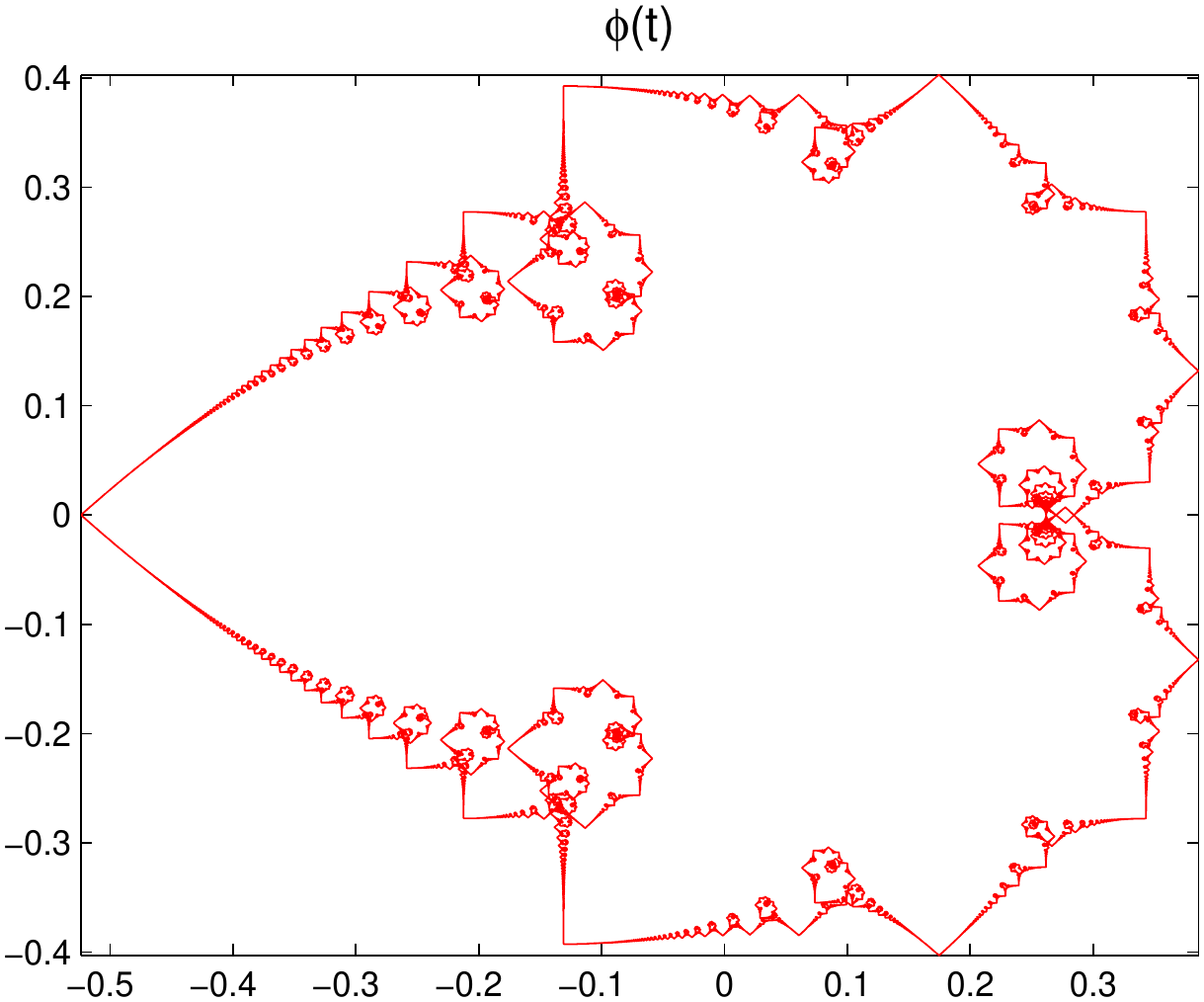}~
\includegraphics[width=0.3333\textwidth, clip=true]{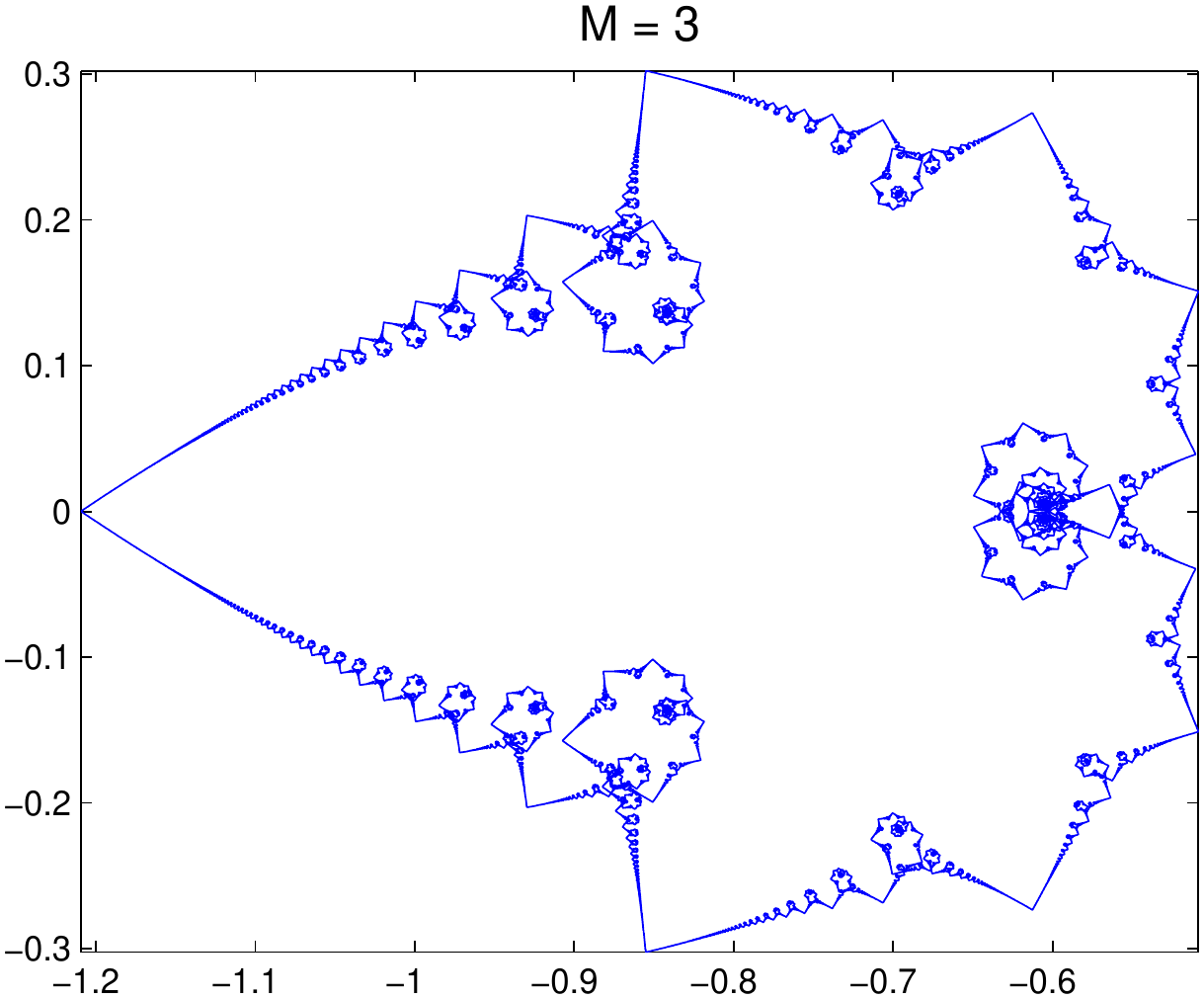}~
\includegraphics[width=0.3333\textwidth, clip=true]{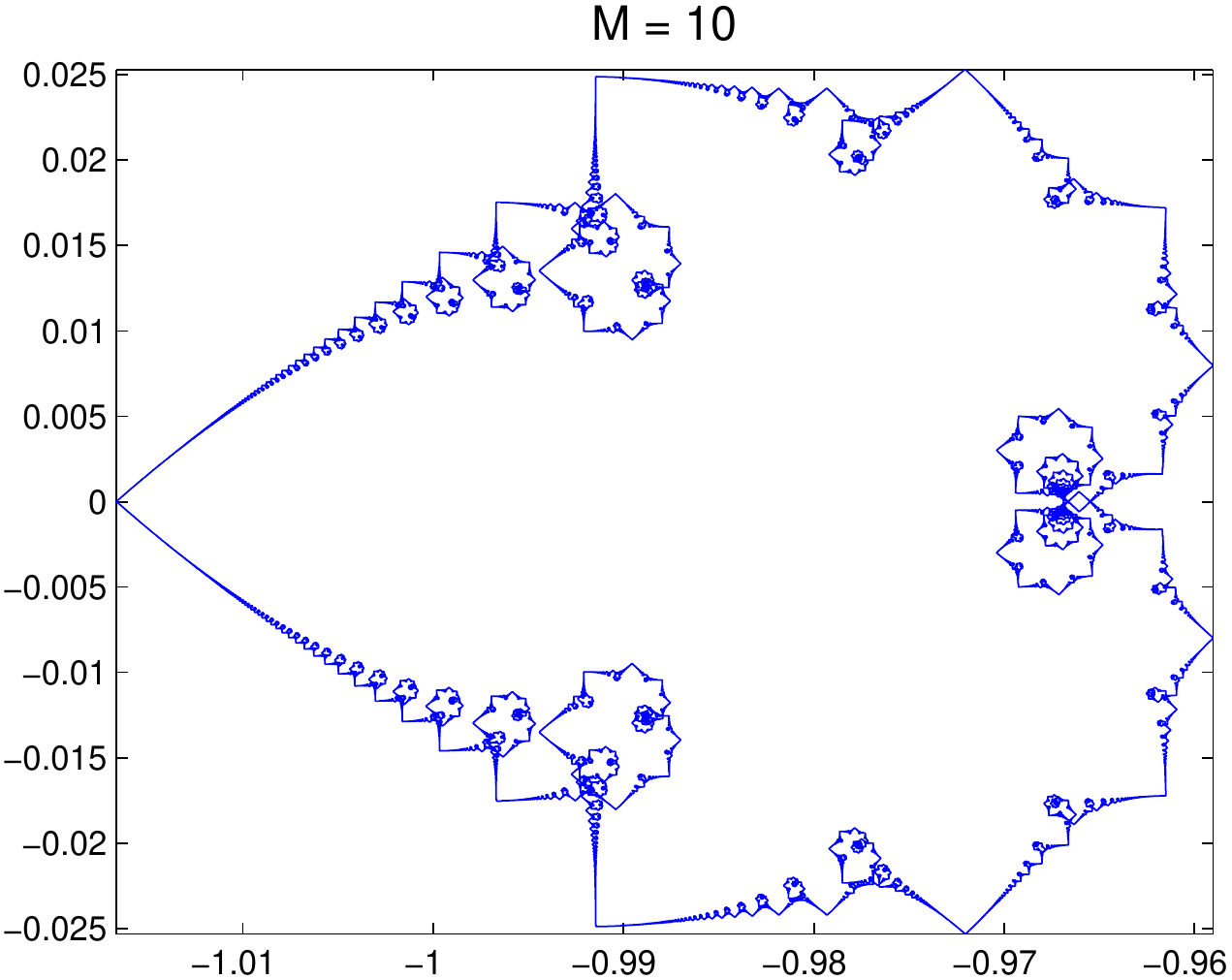}

\caption{$\phi(t)$ (red), as defined in \eref{e:phit2}, against $z_M(t)$ (blue), for $M = 3, 10$. To compute $z_M(t)$, we took $N/M = 8192$; to compute $\phi(t)$, we took $\sum_{k = 0}^{8192-1}$. We have plotted the $z_M(t)$ corresponding to all the $N_t + 1$ points $\X^{(n)}(0)$ obtained during the execution of \eref{e:rungekutta}. We have computed algebraically and plotted $\phi(t)$ at those times. Observe that, except for the scaling, $\phi(t)$ and $z_{10}(t)$ are visually undistinguishable.}\label{f:riemann}
\end{figure}

\noindent In Figure \ref{f:Xs0}, we have plotted the numerically obtained $z(t)$ for $M = 3, 4, 5$. Besides the conspicuous fractal character of the curves, which immediately reminds us of the works \cite{Peskin94,Stern94}, their most striking feature (see Figure \ref{f:riemann}) is how much $z(t) - \rmi c_M t$, i.e., z(t) without the vertical movement, resembles the graph of
\begin{equation}
\label{e:phit1}
\phi(t) = \sum_{k = 1}^{\infty}\frac{\rme^{\pi \rmi k^2 t}}{\rmi\pi k^2}, \quad t\in[0,2].
\end{equation}

\noindent The function $\phi(t)$ was conveniently used in \cite{Du}, in order to study Riemann's non-differen\-tiable function, which is precisely the real part of $\phi(t)$, i.e.,
\begin{equation}
\label{e:f(t)}
f(t) = \sum_{k = 1}^{\infty}\frac{\sin(\pi k^2 t)}{\pi k^2}.
\end{equation}

\noindent This function is non-differentiable everywhere, except at those rational points $t = p / q$, with $p$ and $q$ both odd. On the other hand, when constructing algebraically $\X(0, t)$ for a given t, we have observed that the only times at which $\X(0, t)$ has no corner are of the form $t_{pq} = (2\pi/M^2)(p/q)$, with $p$ odd and $q\equiv 2\bmod 4$. Therefore, in order to compare $\phi(t)$ and $\X(0, t)$, we have to redefine slightly $\phi(t)$ in \eref{e:phit1}:
\begin{equation}
\label{e:phit2}
\phi(t) = -\sum_{k = 1}^{\infty}\frac{\rme^{-2\pi \rmi k^2 t}}{\pi k^2}, \quad t\in[0,1],
\end{equation}

\noindent i.e., we multiply $\phi$ in \eref{e:phit1} by $-\rmi $, in order to orientate it correctly; then we change the period from $t\in[0, 2]$ to $t\in[0, 1]$, and the orientation of the parameterization, by doing $t \equiv 2-2t$. On the other hand, given a number of initial sides $M$, we define
\begin{equation}
z_M(t) \equiv z(\tfrac{2\pi t}{M^2}) - \rmi c_M\tfrac{2\pi t}{M^2}, \quad t\in[0, 1].
\end{equation}

\noindent In Figure \ref{f:riemann}, we have plotted $\phi(t)$ (red), as defined in \eref{e:phit2}, versus $z_M(t)$, for $M = 3, 10$. Since we have obtained numerically $z_M(t)$ at $N_t+1$ points in $[0, 1]$, we have evaluated $\phi(t)$ precisely at those points. Although the figures are not identical, they are extremely similar. Moreover, if we scale them accordingly, we observe that, the bigger $M$ grows, the more similar to $\phi(t)$ is $z_M(t)$. Indeed, $\phi(t)$ and $z_{10}(t)$, except for the scaling, are visually undistinguishable.

\begin{figure}[t!]
\center

\includegraphics[width=0.5\textwidth, clip=true]{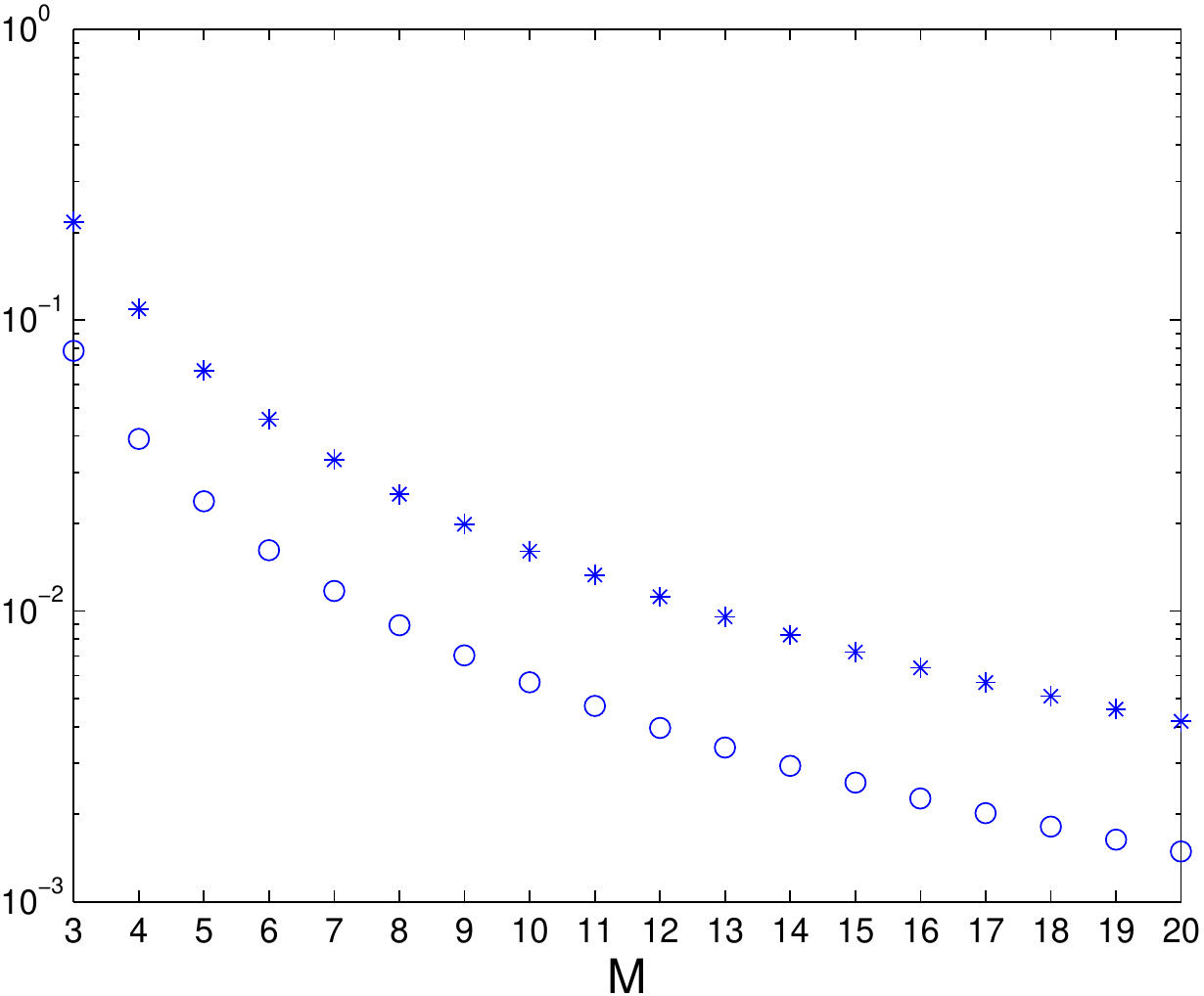}

\caption{$\max_t|\phi(t) - \lambda_Mz_M(t) - \mu_M|$ (starred line) and $\max_t|(\phi(t) - \lambda_Mz_M(t) - \mu_M) / \phi(t)|$ (circled line), with $\lambda_M$ and $\mu_M$ given by \eref{e:lambdamu}. Bearing in mind that we are comparing $38707200+1$ points one by one, the convergence is pretty strong.} \label{f:convergence}
\end{figure}

In order to give some numerical evidence that an adequately scaled $z_M(t)$ converges to $\phi(t)$, we study $\phi(t) - \lambda_Mz_M(t) - \mu_M$, for certain $\lambda_M\in\mathbb R$ and $\mu_M\in\mathbb C$. After applying a least-square fitting, we choose
\begin{equation}
\label{e:lambdamu}
\eqalign{
\lambda_M = \Re\left(\frac{\mean\left[(z_M(t) - \mean(z_M(t)))(\bar\phi(t) - \mean(\bar\phi(t)))\right]}{\mean(|z_M(t) - \mean(z_M(t))|^2)}\right),
    \\
\mu_M = \mean(\phi(t)) - \lambda_M\mean(z_M(t)),
}
\end{equation}

\noindent where $\mean(z_M(t))$ denotes the mean over the $N_t$ + 1 points, etc. Figure \ref{f:convergence} shows the absolute and relative discrepancies between $\phi(t)$ and the scaled $z_M(t)$, i.e., $\max_t|\phi(t) - \lambda_Mz_M(t) - \mu_M|$ and $\max_t|(\phi(t) - \lambda_Mz_M(t) - \mu_M) / \phi(t)|$, for $M = 3, 4, \ldots, 20$. Bearing in mind that we are comparing $38707200+1$ different points one by one, the results give a pretty strong evidence of convergence as the number of initial sides $M$ tends to infinity, and strongly support the idea that $z_M$ is also a multifractal, for all $M \ge 3$.

As we mentioned in the introduction, the proof that $f(t)$, as defined in \eref{e:f(t)}, is a multifractal was done in \cite{Ja}. In this respect, it is fundamental the identity given for $\phi$, as defined in \eref{e:phit1}, which states that
\begin{equation}
\label{e:phiduistermaat}
\phi(t) = \phi(t_{pq}) + \rme^{\pi \rmi m / 4}q^{-1/2}(t - t_{pq})^{1/2} + \mbox{ lower-order terms},
\end{equation}

\noindent with $p, q\in\mathbb Z$, $q > 0$, $\gcd(p, q) = 1$, $m\equiv m(t_{pq})\in\mathbb Z/8\mathbb Z$. This identity is proved in Theorem 4.2 in \cite{Du}. Of particular relevance is the fact that the H\"older exponent is $1/2$, i.e., that $|\phi(t) - \phi(t_{pq})| = q^{-1/2}|t - t_{pq}|^{1/2} + $ lower-order terms.

\begin{figure}[t!]
\center\includegraphics[width=0.5\textwidth, clip=true]{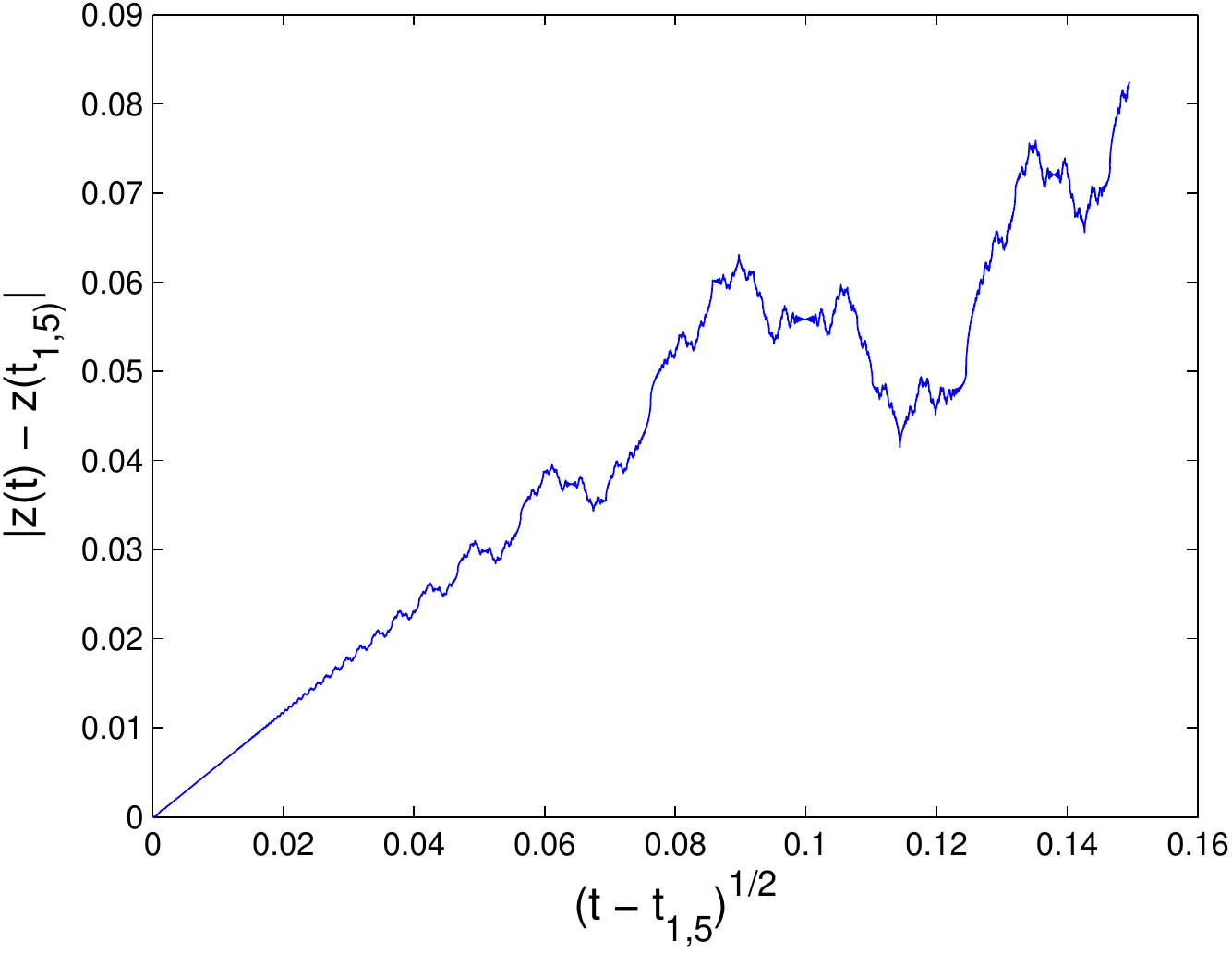}
\caption{$|z(t) - z(t_{1,5})|$ versus $(t - t_{1,5})^{1/2}$, for $M = 3$. The asymptotically linear relation between both quantities as $t\to t_{1,5}$ is evident.}\label{f:holder}
\end{figure}

In order to prove analytically that $z_M(t)$, is a multifractal, we would need an equivalent of \eref{e:phiduistermaat}. We have made some numerical experiments (see for example Figure \ref{f:holder}, for $M = 3$, $p = 1$, $q = 5$), and all of them give strong evidence that the H\"older exponent of $z(t)$ is $1/2$ for rational times, i.e., that $|z_M(t) - z_M(t_{pq})| = \mathcal O(|t - t_{pq}|^{1/2})$.

Nevertheless, how the constants depend on the denominator $q$, something which is a fundamental ingredient in the arguments in \cite{Du} and in \cite{Ja}, is unclear. This question deserves a much more detailed analysis that we plan to make in a forthcoming paper.

\subsection{$\T(s, t_{pq})$, for $q\gg1$}

\label{s:fractalsT}

In the previous subsection, we have given evidence of the multifractal character of the curve $\X(0, t)$. However, in \cite{DelahozGarciaCerveraVega09}, fractal-like phenomena were observed also on the tangent vector $\T$ when imposing fixed boundary conditions. In the following lines, we will show that similar observations are valid for periodic boundary conditions, too.

As we have seen, $\X$ can be recovered algebraically up to a vertical movement, while $\T$ can be completely recovered algebraically at rational times. Furthermore, $\T_{alg}$ is really a piecewise constant function, with no noise associated to the Gibbs phenomenon. This gives us a very powerful tool to progress in the understanding of $\T$ (and also of $\X$), avoiding numerical simulations at all.

A very interesting question is, given a rational time $t_{pq}$, with \emph{small} $q$, corresponding in $\X$ to a polygon with $Mq$ or $Mq/2$ sides, and in $\T$ to a piecewise continuous function with $Mq$ or $Mq/2$ jumps, what happens at a time $t_{pq} + \varepsilon$, $|\varepsilon|\ll1$? Let us take $\varepsilon = (2\pi/M^2)/q'$, with $q'\gg1$, in order that $|\varepsilon| \ll 1$. Assuming that $\gcd(q, q') = 1$, then, $\tfrac{p}{q} + \tfrac{1}{q'} = \tfrac{pq'+1}{qq'}$, so $t_{pq} + \varepsilon$ would correspond to a polygon with $qq'$ or $qq'/2$ sides.

\begin{figure}[t!]
\center
\includegraphics[width=0.5\textwidth, clip=true]{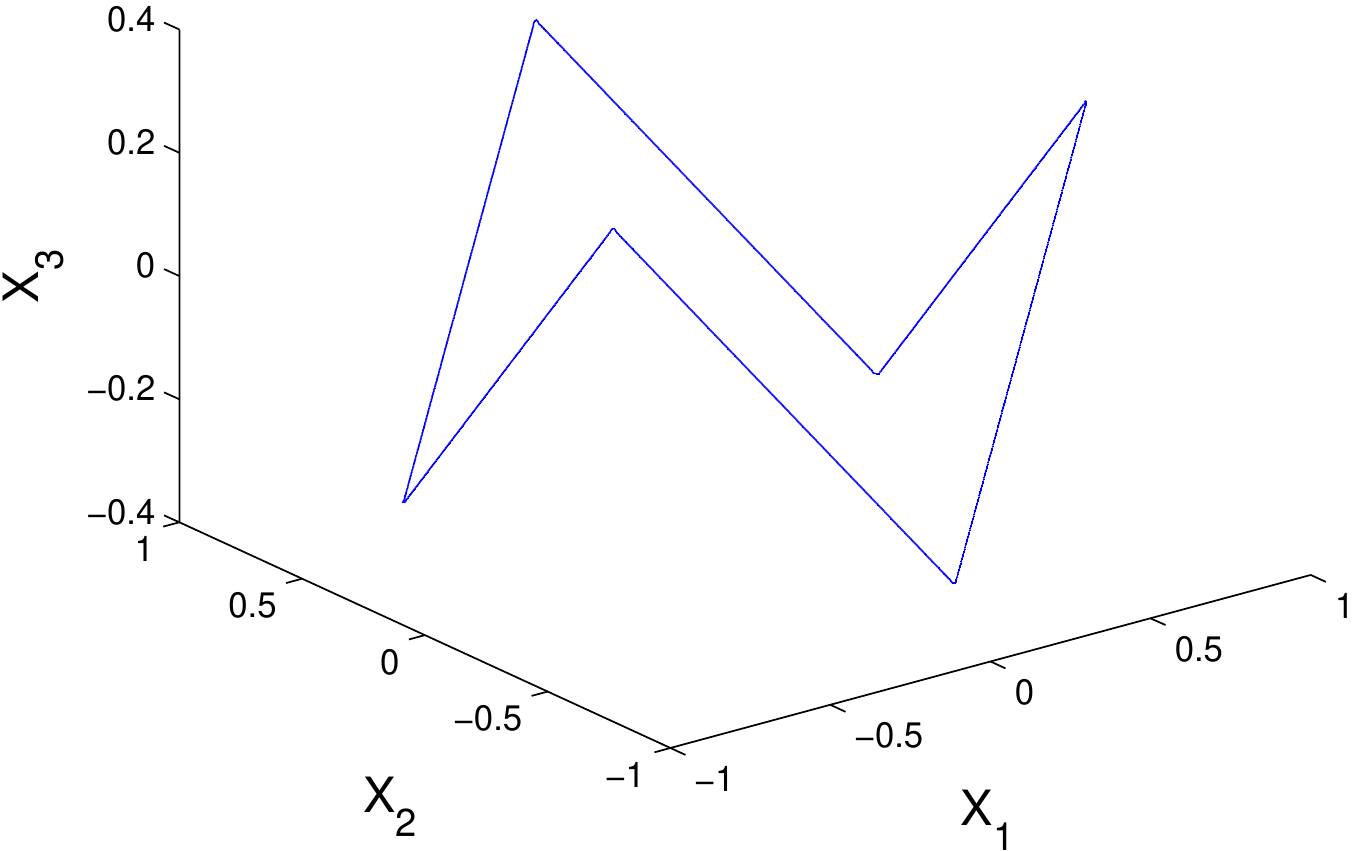}~
\includegraphics[width=0.5\textwidth, clip=true]{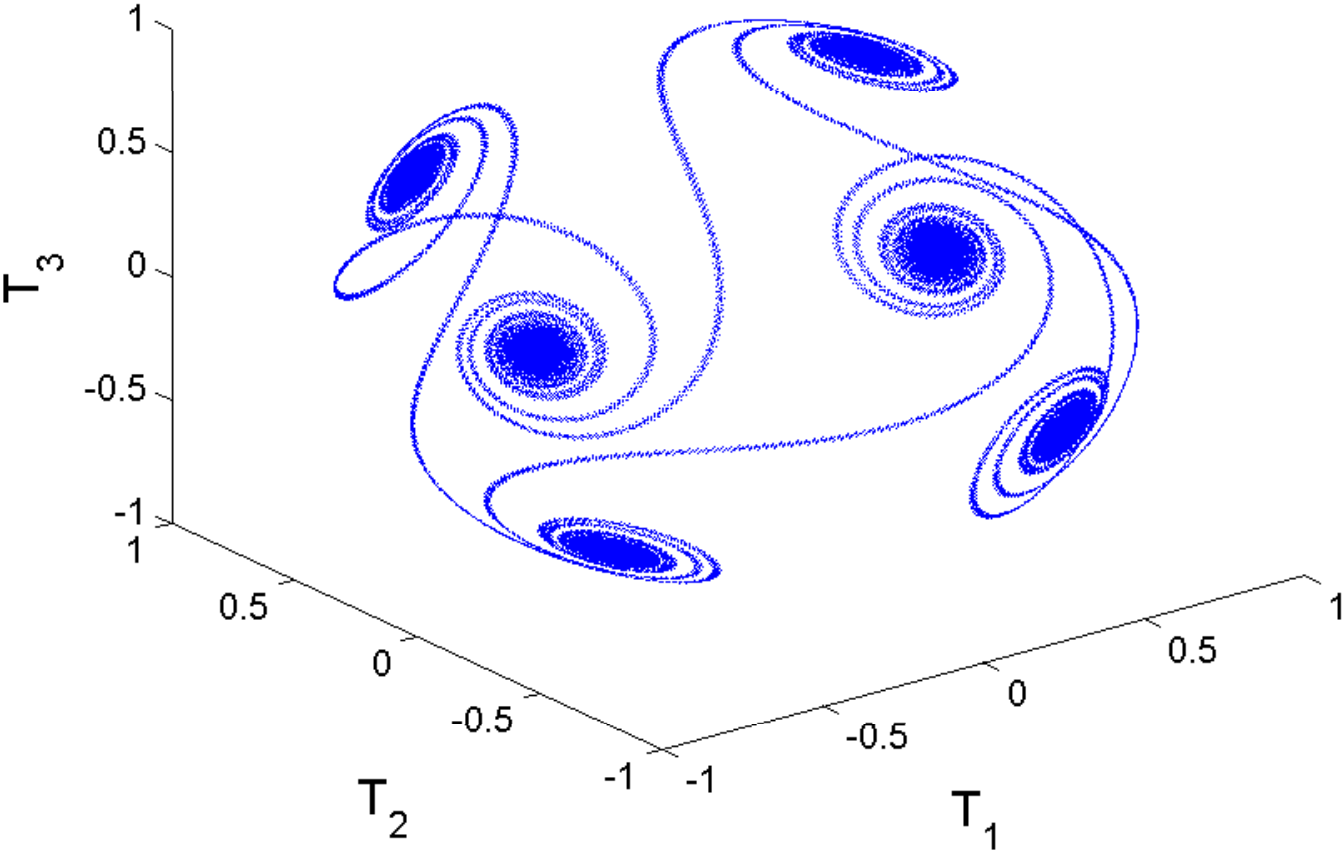}
\caption{$\X_{alg}$ and $\T_{alg}$, at $t = \tfrac{2\pi}{9}(\tfrac{1}{4}+\tfrac{1}{49999}).$}\label{f:XTeps}
\end{figure}

In Figure \ref{f:XTeps}, we have have plotted $\X_{alg}$ and $\T_{alg}$, for $M = 3$, at $t = \tfrac{2\pi}{9}(\tfrac{1}{4}+\tfrac{1}{49999}) = \tfrac{2\pi}{9}\cdot\tfrac{50003}{199996}$. While visually there is no difference whatsoever between $\X_{alg}$ at $t_{1,4} = \tfrac{2\pi}{9}\cdot\tfrac{1}{4}$ and $\X_{alg}$ at $t = \tfrac{2\pi}{9}\cdot\tfrac{50003}{199996}$, we do not have a skew polygon with 6 sides, but a skew polygon with $Mq/2 = 299994$ sides that closely resembles a polygon with 6 sides. On the other hand, the corresponding $\T_{alg}$ are very different. Indeed, in the plot of $\T_{alg}$, we can clearly appreciate six spiral-like structures whose centers are precisely the six constant values of $\T$ at $t_{1,4}$. Nevertheless, it is important to underline that these structures are not really spirals, but the plot of the $299994$ different values taking by $\T$, that closely resembles a curve with six spirals. Furthermore, these spirals remind us of the Cornu spirals that appeared in \cite{gutierrez}. An open question that arises naturally is up to what extent the multiple-corner problem can be explained as a sum of several one-corner problems.

\begin{figure}[t!]
\center
\includegraphics[width=0.5\textwidth, clip=true]{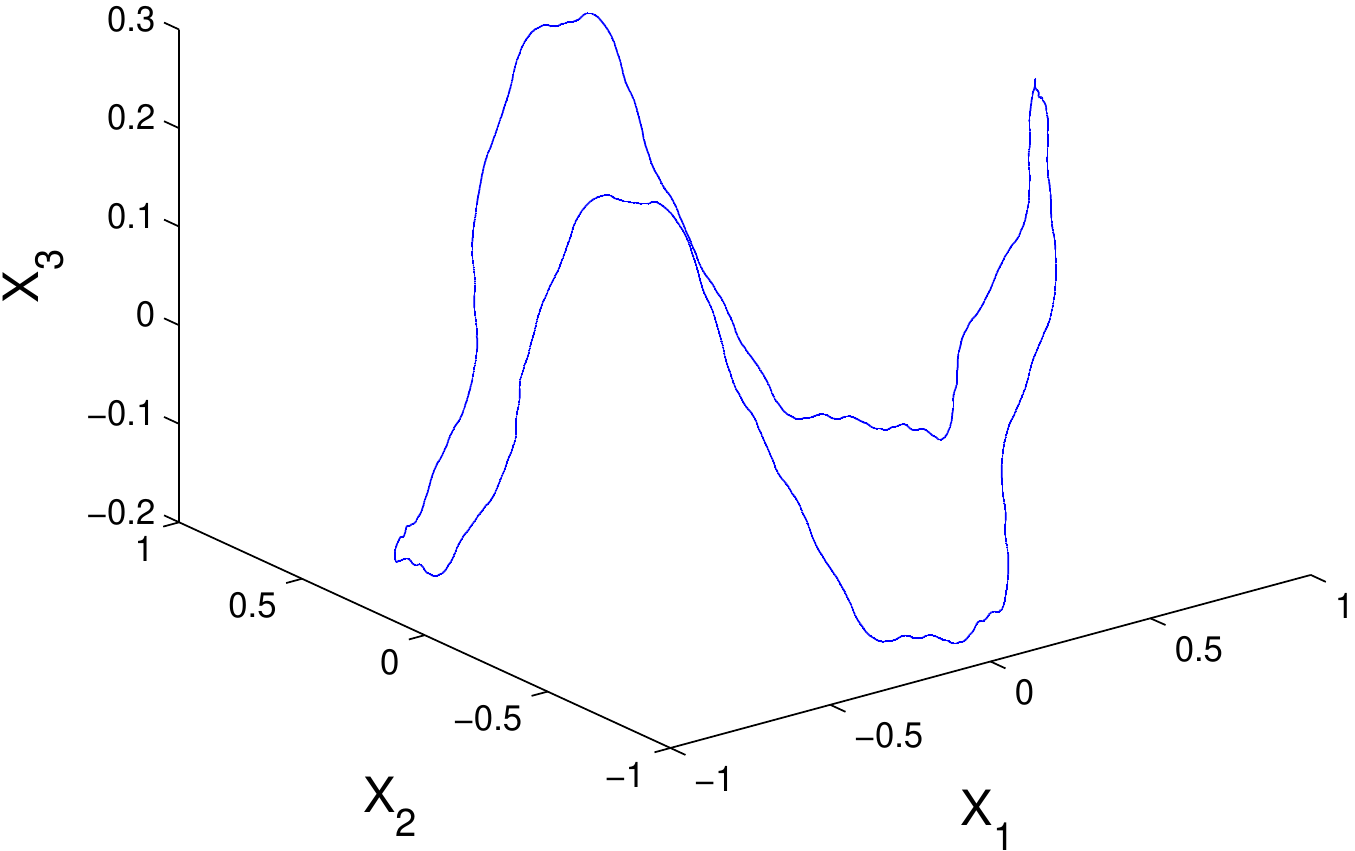}~
\includegraphics[width=0.5\textwidth, clip=true]{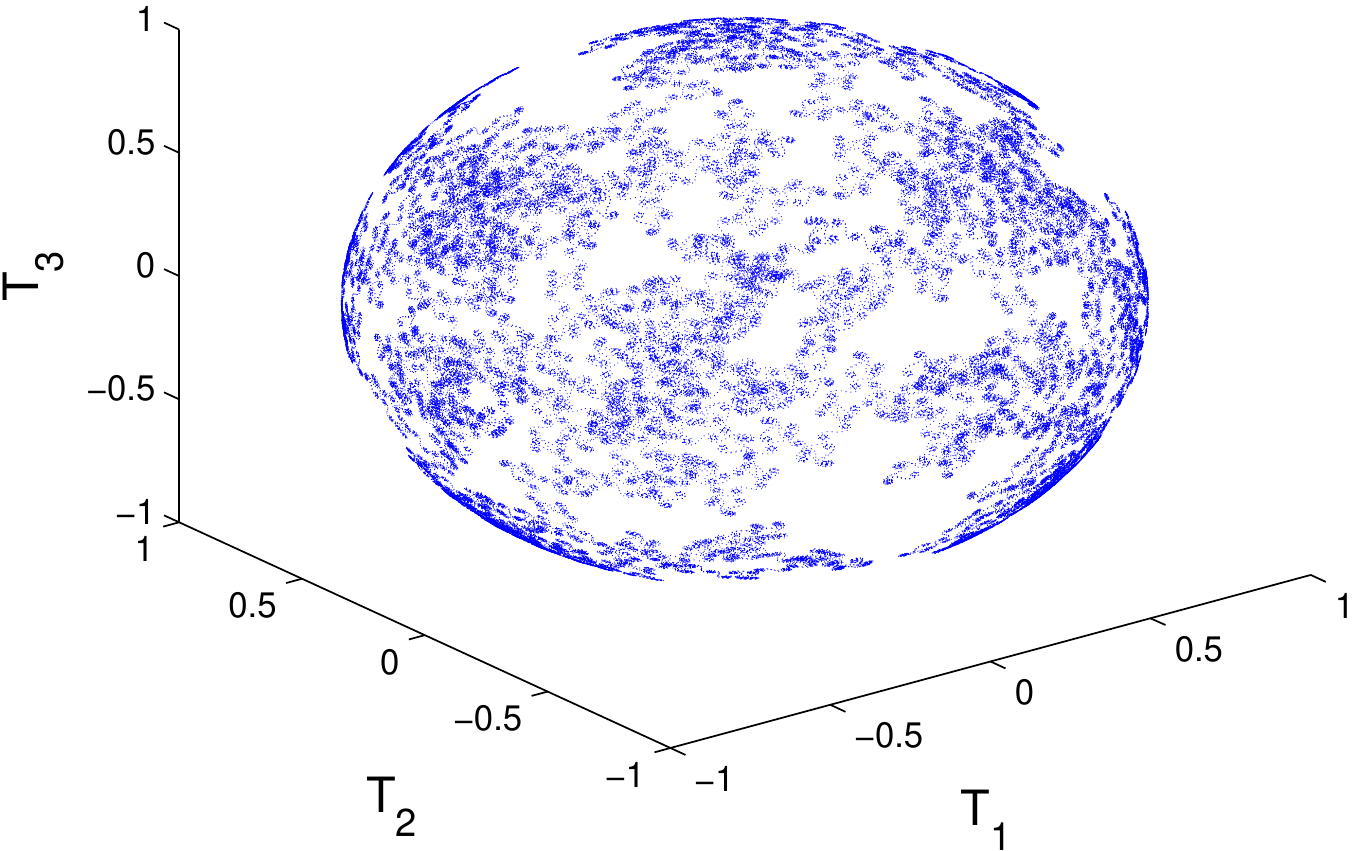}
\caption{$\X_{alg}$ and $\T_{alg}$, at $t = \tfrac{2\pi}{9}(\tfrac{1}{4} + \tfrac{1}{41} + \tfrac{1}{401}) = \tfrac{2\pi}{9}\cdot\tfrac{18209}{65764}$. }\label{f:XTeps2}
\end{figure}

\begin{figure}[t!]
\center
\includegraphics[width=0.5\textwidth, clip=true]{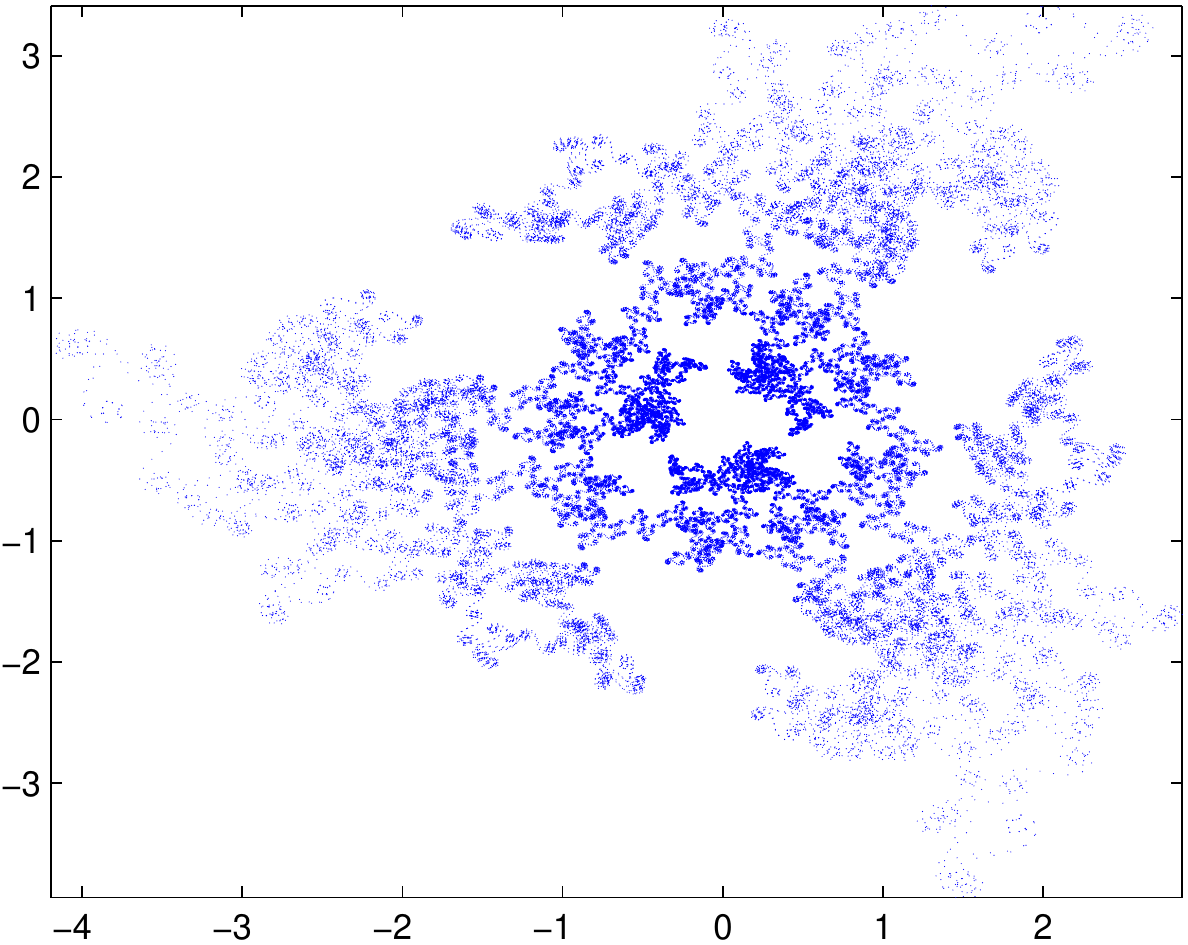}
\caption{Stereographic projection of the right-hand side of Figure \ref{f:XTeps2}.}\label{f:Tepszoom}
\end{figure}

Another interesting question is what happens if we take a rational time $t_{pq}$ with \emph{large} $q$, such that there is no pair $(p', q')$, with both $q'$ and $|\tfrac{p}{q} - \tfrac{p'}{q'}|$ \emph{small}. In this case, the situation is very different. In Figure \ref{f:XTeps2}, we have plotted $\X_{alg}$ and $\T_{alg}$ for $M = 3$, at $t = \tfrac{2\pi}{9}(\tfrac{1}{4} + \tfrac{1}{41} + \tfrac{1}{401}) = \tfrac{2\pi}{9}\cdot\tfrac{18209}{65764}$. While the left-hand side is not so different from the left-hand side of Figure \ref{f:XTeps}, the right-hand side is a set of $Mq /2 = 98646$ points in $\mathbb S^2$ that creates a spectacular fractality sensation with spiral-like structures at three or four different scales. This can be better appreciated in Figure \ref{f:Tepszoom}, where we have plotted the stereographic projection \eref{e:stereo} of $\T$ onto $\mathbb C$.

In general, a rational time $t_{pq}$, with $q\gg1$, can be regarded as an approximation of an irrational time.
Trying to fully understand the phenomena exhibited by $\T(s, t)$ at those times is a challenging question that we postpone for the future.

\section{Conclusions}

\label{s:conclusions}

In this paper, we have studied the evolution of \eref{e:xt} and \eref{e:schmap}, for a regular planar polygonal with $M$ sides as the initial datum. The algebraic calculations, backed by complete numerical simulations, prove (under a certain uniqueness assumption) that the curve $\X(s, t)$ is a polygon at times which are rational multiples of $2\pi/M^2$, i.e., $t_{pq}=(2\pi/M^2)(p/q)$, $\gcd(p, q) = 1$, with the number of sides depending on $q$, while the tangent vector $\T(s, t)$ is piecewise constant at those times. This could be seen as a nonlinear version of the Talbot effect, in the spirit of the work \cite{BG} by Berry and Goldberg.

We have obtained a striking connection between the curve $\X(0, t)$ and the so-called Riemann's non-differentiable function. In \cite{Ja}, S. Jaffard proved that this function is an example of multifractal whose spectrum of singularities satisfies the Frisch-Parisi conjecture. Although there is strong numerical evidence that $\X(0, t)$ is also a multifractal, an analytical proof seems to be challenging. In fact, a first step in this direction is to give sense to our solutions from an analytical point of view.

In the future, we also plan to extend these ideas to arbitrary polygons, and to do a more detailed study on fractality; for such purposes, a more complete algebraically description of $\X$ and $\T$ is no doubt required. Finally, it would be interesting to prove \eref{e:cosrho}.

\section*{Acknowledgements}

We are indebted to the referees and editors for the thorough, constructive and helpful comments and suggestions. We also wish to thank V. Banica and F. Chamizo for enlightening conversations.

This work was supported by the Spanish Ministry of Economy and Competiveness, with the project MTM2011-24054, and by the Basque Government, with the project IT641-13.

\appendix

\section{Generalized Quadratic Gau{\ss} Sums}

\label{s:gausssums}

The generalized quadratic Gau{\ss} sums are defined by
\begin{equation}
G(a, b, c) = \sum_{l=0}^{|c| - 1}\rme^{2\pi \rmi (al^2 + bl)/c},
\end{equation}

\noindent for given integers $a, b, c$, with $c\not=0$. From now on, we assume $c > 0$, and $\gcd(a, c) = 1$, which are the cases dealt with in this paper. The value of these sums was calculated for the first time by Gau{\ss} for the case $b = 0$ \cite{Berndt1988}. More precisely, given two integers $a, c$ such that $\gcd(a, c) = 1$, then
\begin{equation}
G(a, 0, c) = \sum_{l=0}^{c - 1}\rme^{2\pi \rmi al^2/c} =
\cases{
(\tfrac{c}{a})(1+\rmi ^a)\sqrt c,  & if $c \equiv 0 \bmod 4$,
    \\
(\tfrac{a}{c})\sqrt c, & if $c \equiv 1 \bmod 4$,
    \\
0, & if $c \equiv 2 \bmod 4$,
    \\
(\tfrac{a}{c})\rmi\sqrt c, & if $c \equiv 3 \bmod 4$,
}
\end{equation}

\noindent where $(\tfrac{a}{c})$ denotes the Jacobi symbol. As we will see in the following lines, the generalized quadratic Gau{\ss} sums can be reduced to normal quadratic Gau{\ss} sums by completing the square. Since they are multiplicative, i.e.,
\begin{equation}
G(a, b, cd) = G(ac, b, d)G(ad, b, c), \quad \mbox{with } \gcd(c, d) = 1,
\end{equation}

\noindent we can assume without loss of generality that $c$ is either odd or a power of two. If $c$ is odd, let us find a certain $\phi(a)$ such that $4a\phi(a)\equiv 1 \bmod c$. This $\phi(a)$ is unique modulo $c$ and is precisely the inverse of $4a$ in $\mathbb Z/ c\mathbb Z$; since $\gcd(4a, c) = 1$, its existence is guaranteed by B\'ezout's lemma, and it can be efficiently computed, for instance, by the extended Euclidean algorithm. Then, we have
\begin{eqnarray}
G(a, b, c) & = \sum_{l=0}^{c - 1}\rme^{2\pi \rmi (al^2 + 4a\phi(a)bl + (4a\phi(a))\phi(a)b^2 - \phi(a)b^2)/c}
    \cr
& = \rme^{-2\pi \rmi \phi(a)b^2/c}\sum_{l=0}^{c - 1}\rme^{2\pi \rmi  al^2/c}
    \cr
& =
\cases{
\rme^{-2\pi \rmi \phi(a)b^2/c}(\tfrac{a}{c})\sqrt c, & if $c \equiv 1 \bmod 4$,
    \\
\rme^{-2\pi \rmi \phi(a)b^2/c}(\tfrac{a}{c})\rmi\sqrt c, & if $c \equiv 3 \bmod 4$.
}
\end{eqnarray}

\noindent If $c$ is a power of two, we first observe that
\begin{eqnarray*}
G(a, b, c) & = \sum_{l = 0}^{c/2 - 1} \rme^{2\pi \rmi (al^2 + bl) / c} + \sum_{l = 0}^{c/2 - 1} \rme^{2\pi \rmi (a(l + c/2)^2 + b(l + c/2)) / c}
    \\
& =
\cases{
0, & if $a(c/2) \not\equiv b\bmod 2$,
    \\
2\sum_{l = 0}^{c/2 - 1} \rme^{2\pi \rmi (al^2 + bl) / c}, & if $a(c/2) \equiv b\bmod 2$;
}
\end{eqnarray*}

\noindent therefore, we have two cases: $c = 2$ and $c > 2$. If $c = 2$, then $a$ is odd, so
\begin{equation}
G(a, b, 2) =
\cases{
0, & if $b$ even,
    \\
2, & if $b$ odd.
}
\end{equation}

\noindent On the other hand, if $c > 2$, $G(a, b, c)=0$, if $b$ is odd. Let us suppose $b$ is even. Then, since $a$ is odd, we follow the previous reasoning to find a certain $\phi(a)$ such that $a\phi(a)\equiv 1\bmod c$, so
\begin{eqnarray}
G(a, b, c) & = \sum_{l=0}^{c - 1}\rme^{2\pi \rmi (al^2 + a\phi(a)2(b/2)l + (a\phi(a))\phi(a)(b/2)^2 - \phi(a)(b/2)^2)/c}
    \cr
& = \rme^{-\pi \rmi \phi(a)b^2/(2c)}\sum_{l=0}^{c - 1}\rme^{2\pi \rmi  al^2/c}
    \cr
& = \rme^{-\pi \rmi \phi(a)b^2/(2c)}(\tfrac{c}{a})(1+\rmi ^a)\sqrt c.
\end{eqnarray}

\noindent In general, given an arbitrary $c\in\mathbb N$, we can factorize it as $c = 2^rc'$ and use the multiplicative character of the generalized quadratic Gau{\ss} sums to calculate $G(a, b, c)$. For instance, if $c$ is even, but $c / 2$ is odd, then
\begin{eqnarray}
G(a, b, c) & = G(a, b, 2(c/2))
    \cr
& = G(a(c/2), b, 2)G(2a, b, c/2)
    \cr
& = \cases{
0, & if $a(c/2) \not\equiv b\bmod 2$,
    \\
2G(2a, b, c/2), & if $a(c/2) \equiv b\bmod 2$,
}
\end{eqnarray}

\noindent i.e, $G(a, b, c) = 0$, if $b$ is even, etc.

To conclude this appendix, let us mention explicitly the value of $|G(a, b, c)|$, deduced from the previous calculations, and which is of especial relevance in this paper:
\begin{equation}
\label{e:absGabc}
|G(a, b, c)| = \cases{
\sqrt c, & if $c$ odd,
    \\
\sqrt{2c}, & if $c$ even $\wedge$ $c/2\equiv b\bmod 2$,
    \\
0, & if $c$ even $\wedge$ $c/2\not\equiv b\bmod 2$.
}
\end{equation}


\section*{References}

\begin{thebibliography}{10}

\bibitem{darios}
L.~S.~Da Rios.
\newblock {Sul moto d'un liquido indefinito con un filetto vorticoso di forma
  qualunque}.
\newblock {\em Rend. Circ. Mat. Palermo}, 22(1):117--135, 1906.
\newblock In Italian.

\bibitem{arms}
R.~J. Arms and F.~R. Hama.
\newblock {Localized-Induction Concept on a Curved Vortex and Motion of an
  Elliptic Vortex Ring}.
\newblock {\em Phys. Fluids}, 8(4):553--559, 1965.

\bibitem{batchelor}
G.~K. Batchelor.
\newblock {\em {An Introduction to Fluid Dynamics}}.
\newblock Cambridge Mathematical Library. Cambridge University Press, 1967.

\bibitem{saffman}
P.~G. Saffman.
\newblock {\em {Vortex Dynamics}}.
\newblock Cambridge Monographs on Mechanics. Cambridge University Press, 1995.

\bibitem{landau}
L.~D. Landau and E.~M. Lifshitz.
\newblock {On the theory of the dispersion of magnetic permeability in
  ferromagnetic bodies}.
\newblock {\em Phys. Z. Sowjet.}, 8(2):153--169, 1935.

\bibitem{delahoz2007}
F.~de~la Hoz.
\newblock {Self-similar solutions for the 1-D Schr\"odinger map on the
  hyperbolic plane}.
\newblock {\em Math. Z.}, 257(1):61--80, 2007.

\bibitem{hasimoto}
H.~Hasimoto.
\newblock {A soliton on a vortex filament}.
\newblock {\em J. Fluid Mech.}, 51(3):477--485, 1972.

\bibitem{gutierrez}
S.~Guti\'errez, J.~Rivas, and L.~Vega.
\newblock {Formation of singularities and self-similar vortex motion under the
  localized induction approximation}.
\newblock {\em Comm. PDE}, 28(5--6):927--968, 2003.

\bibitem{buttke87}
T.~F. Buttke.
\newblock {A Numerical Study of Superfluid Turbulence in the Self-Induction
  Approximation}.
\newblock {\em J. Comput. Phys.}, 76(2):301--326, 1998.

\bibitem{DelahozGarciaCerveraVega09}
F.~de~la Hoz, C.~J. {Garc\'\i a-Cervera}, and L.~Vega.
\newblock {A Numerical Study of the Self-Similar Solutions of the Schr\"odinger
  Map}.
\newblock {\em SIAM J. Appl. Math.}, 70(4):1047--1077, 2009.

\bibitem{Peskin94}
C.~S. Peskin and D.~M. McQueen.
\newblock Mechanical equilibrium determines the fractal fiber architecture of
  the aortic heart valve leaflets.
\newblock {\em Am. J. Physiol.}, 266(1):H319--H328, 1994.

\bibitem{Stern94}
J.~V. Stern and C.~S. Peskin.
\newblock Fractal dimension of an aortic heart valve leaflet.
\newblock {\em Fractals}, 2(3):461--464, 1994.

\bibitem{BG}
M.~V. Berry and J.~Goldberg.
\newblock {Renormalisation of curlicues}.
\newblock {\em Nonlinearity}, 1(1):1--26, 1988.

\bibitem{BK}
M.~V. Berry and S.~Klein.
\newblock {Integer, fractional and fractal Talbot effects}.
\newblock {\em J. Mod. Optics}, 43:2139--2164, 1996.

\bibitem{BMS}
M.~V. Berry, I.~Marzoli, and W.~Schleich.
\newblock {Quantum carpets, carpets of light}.
\newblock {\em Physics World}, 14(6):39--44, 2001.

\bibitem{B1}
M.~V. Berry.
\newblock {Quantum fractals in boxes}.
\newblock {\em J. Phys. A: Math. Gen.}, 29:6617--6629, 1996.

\bibitem{ET}
M.~B. Erdo{\u g}an and N.~Tzirakis.
\newblock {Talbot effect for the cubic nonlinear Schr\"odinger equation on the
  torus}.
\newblock 2013.
\newblock arXiv:1303.3604 [math.AP].

\bibitem{Ol}
P.~J. Olver.
\newblock {Dispersive quantization}.
\newblock {\em Amer. Math. Monthly}, 117(7):599--610, 2010.

\bibitem{Os}
K.~I. Oskolkov.
\newblock {A Class of I.M. Vinogradov's Series and Its Applications in Harmonic
  Analysis}.
\newblock {\em Progress in Approximation Theory}, 19:353--402, 1992.

\bibitem{KR}
L.~Kapitanski and I.~Rodnianski.
\newblock {Does a Quantum Particle Know the Time?}
\newblock {\em Emerging Applications of Number Theory, IMA Vol. Math. Appl.},
  109:355--371, 1999.

\bibitem{CET}
V.~Chousionis, M.~B. Erdo{\u g}an, and N.~Tzirakis.
\newblock {Fractal solutions of linear and nonlinear dispersive partial
  differential equations}.
\newblock 2014.
\newblock arXiv:1406.3283v1 [math.AP].

\bibitem{BV1}
V.~Banica and L.~Vega.
\newblock {Scattering for 1D cubic NLS and singular vortex dynamics}.
\newblock {\em J. Eur. Math. Soc. (JEMS)}, 14(1):209--253, 2012.

\bibitem{BV2}
V.~Banica and L.~Vega.
\newblock {Stability of the selfsimilar dynamics of a vortex filament}.
\newblock 2012.
\newblock arXiv:1202.1106 [math.AP].

\bibitem{BV3}
V.~Banica and L.~Vega.
\newblock {The initial value problem for the Binormal Flow with rough data}.
\newblock 2013.
\newblock preprint.

\bibitem{didier}
R.~L. Jerrard and D.~Smets.
\newblock {On the motion of a curve by its binormal curvature}.
\newblock {\em Jour. Eur. Math. Soc.}, 2014.
\newblock To appear.

\bibitem{Didier2}
R.~L. Jerrard and D.~Smets.
\newblock {On Schr\"odinger maps from $T^1$ to $S^2$}.
\newblock {\em Ann. Sci. \'Ec. Norm. Sup\'er. (4)}, 45(4):637--680, 2013.

\bibitem{GG}
F.~F. Grinstein and E.~J. Gutmark.
\newblock {Flow control with noncircular jets}.
\newblock {\em Ann. Rev. Fluid Mech.}, 31:239--272, 1999.

\bibitem{GGP}
F.~F. Grinstein, E.~J. Gutmark, and T.~Parr.
\newblock {Nearfield dynamics of subsonic, free square jets. A computational
  and experimental study}.
\newblock {\em Phys. Fluids}, 7:1483--1497, 1995.

\bibitem{MMG}
R.~S. Miller, C.~K. Madnia, and P.~Givi.
\newblock {Numerical simulation of non-circular jets}.
\newblock {\em Comput. Fluids}, 24(1):1--25, 1995.

\bibitem{KPV}
C.~E. Kenig, G.~Ponce, and L.~Vega.
\newblock {On the ill-posedness of some canonical dispersive equations}.
\newblock {\em Duke Math. J.}, 106(3):617--633, 2001.

\bibitem{Du}
J.~J. Duistermaat.
\newblock {Selfsimilarity of ``Riemann's Nondifferentiable function''}.
\newblock {\em Nieuw Arch. Wisk.}, (4) 9(3):303--337, 1991.

\bibitem{Ja}
S.~Jaffard.
\newblock {The spectrum of singularities of Riemann's function}.
\newblock {\em Rev. Mat. Iberoamericana}, 12(2):441--460, 1996.

\bibitem{CC}
F.~Chamizo and A.~C\'ordoba.
\newblock {Differentiability and dimension of some fractal Fourier series}.
\newblock {\em Adv. Math.}, 142:335--354, 1999.

\bibitem{FP}
U.~Frisch and G.~Parisi.
\newblock {\em {Fully developped turbulence and intermittency}}.
\newblock Proc. Int. Sch. Phys. ``Enrico Fermi'', North-Holland, Amsterdam,
  1985.

\bibitem{F}
U.~Frisch.
\newblock {\em {Turbulence. The Legacy of A. N. Kolmogorov}}.
\newblock Cambridge University Press, 1995.

\bibitem{Ibragimov}
N.~H. Ibragimov.
\newblock {\em {CRC Handbook of Lie Group Analysis of Differential Equations}},
  volume~2.
\newblock CRC Press, 1995.

\bibitem{trenkler2008}
G.~Trenkler and D.~Trenkler.
\newblock {On the product of rotations}.
\newblock {\em Int. J. Math. Educ. Sci. Tech.}, 39(1):94--104, 2008.

\bibitem{FFTW}
M.~Frigo and S.~G. Johnson.
\newblock {The design and implementation of FFTW3}.
\newblock {\em Proc. IEEE}, 93(2):216--231, 2005.

\bibitem{Berndt1988}
B.~C. Berndt, R.~J. Evans, and K.~S. Williams.
\newblock {\em {Gauss and Jacobi Sums}}.
\newblock Canadian Mathematical Society Series of Monographs and Advanced
  Texts. John Wiley and Sons, Inc., 1998.

\end{thebibliography}
\end{document}